\documentclass[preprint,12pt]{elsarticle}

\usepackage{mwe}
\usepackage{xcolor}
 \usepackage{adjustbox}
 \usepackage{soul}
\usepackage[normalem]{ulem}
\usepackage{multicol}
\usepackage[colorlinks]{hyperref}
\usepackage{graphicx,amsmath,amssymb,amsthm,marvosym,tikz}
\usepackage{tikz}
\usetikzlibrary{backgrounds}
\usetikzlibrary{patterns,shapes.misc, positioning}
\newtheorem{theorem}{Theorem}[section]
\newtheorem{lemma}[theorem]{Lemma}
\newtheorem{corollary}[theorem]{Corollary}

\newtheorem{definition}[theorem]{Definition}
\newtheorem{example}[theorem]{Example}
\newtheorem{observation}[theorem]{Observation}
\newtheorem{proposition}[theorem]{Proposition}
\newtheorem{note}[theorem]{Note}
\allowdisplaybreaks
 \usepackage{comment}
 \usepackage{theoremref}
 \usepackage{hyperref}

\usepackage{bibentry}

\journal{}

\begin{document}
	
\begin{frontmatter}
		
\title{Group vertex magicness of product graphs and trees}
\author[a]{Karthik S}
\ead{karthik_p210136ma@nitc.ac.in}	
  \author[b]{M Sabeel K}
\ead{sabeel.math@gmail.com}
\author[a]{K. Paramasivam}
\ead{sivam@nitc.ac.in}
\address[a]{Department of Mathematics, National Institute of Technology Calicut \\ Kozhikode~\textnormal{673601}, India. }
%\author[d]{A V Prajeesh}
%\ead{prajeesh.math@gmail.com}
\address[b]{Department of Mathematics, TKM College of Engineering, Kollam ~\textnormal{691005}, India.}
%\author[a]{N Kamatchi}
%\ead{kamakrishna77@gmail.com}
%\address[a]{Department of Mathematics, Kamaraj College of Engineering and Technology\\ Virudhunagar~\textnormal{625701}, India.}
%\author[c]{\\S Arumugam}
%\ead{s.arumugam.klu@gmail.com}
%\address[c]{National Centre for Advanced Research
%in Discrete Mathematics,
%Kalasalingam University, Krishnankoil~ \textnormal{626126}, India}
\begin{abstract}
  %  Let $\Gamma$ be an Abelian group and $G=(V(G),E(G))$ be a graph. A mapping $\ell$ from $V(G)$ to $\Gamma$ is said to be $\mathcal{A]$-vertex magic if there is an element $\mu \in \Gamma$ such that ${\sum_{x\in N(v)} \ell(x) =\mu}$ for all vertices $v$ of $G$, where $N(v)$ is the set of vertices adjacent to $v$. Then the graph $G$ has an $\Gamma$-vertex magic labeling with magic constant $\mu$. Further, $G$ is group vertex magic if $G$ is $\Gamma$-vertex magic for any Abelian group  $\Gamma$.

\noindent In this article, some necessary conditions of group vertex magicness of graphs with at least one pendant, group vertex magicness of product graphs are proved. Also, a characterization of group vertex magicness of tress of diameter up to $5$, for all infinite Abelian groups with finitely many torsion elements, is obtained.  
\end{abstract}
\begin{keyword}
	Group vertex magic \sep product graphs \sep torsion element  \sep infinite Abelian group   \\
	%%% PACS codes here, in the form: \PACS code \sep code
	%%% MSC codes here, in the form: \MSC code \sep code
	%%% or \MSC[2008] code \sep code (2000 is the default)
\end{keyword}
\end{frontmatter}
\section{Introduction}
 By a graph $G=(V(G),E(G))$, or $G=(V,E)$, we mean a finite undirected graph without loops or multiple edges. For a vertex $v$ in a graph $G$, the set $N_G(v)$ or shortly $N(v)$ is the {\it {open neighborhood}} of $v$, is defined as $N(v)=\{u : uv\in E(G)\}$, and $|N(v)|=\deg(v)$. Two vertices $u$ and $v$ are adjacent or neighbor when they are joined by an edge. We write $``u\sim v"$ to denote $u$ is adjacent with $v$ and  $``u\nsim v"$ to denote $u$ is not adjacent with $v$ in a graph $G$. For more graph theoretic terminology, we refer to  Chartrand et al.\cite{MR3445306} and Hammack and Imrich \cite{MR2817074}.  
 
   A vertex $v$ is a {\textit{pendant vertex}}, $\deg(v)=1$, and the unique vertex adjacent to $v$ is a {\it{support~vertex}}.  A vertex $v$  is a \textit{strong support vertex} if it is adjacent to two or more pendant vertices, and a vertex $v$ is a \textit{weak support vertex} if it is adjacent to a unique pendant.  Also, a vertex $v$ is neither a  support vertex nor a pendant vertex is called a {\it neutral vertex}. Further, $|N(v)| = |N_n(v)|+|N_s(v)|+|N_p(v)|=\deg_n(v)+\deg_s(v)+\deg_p(v)=\deg(v)$, where $N_n(v)$ is the set of all neutral neighbors of $v$, $N_s(v)$ is the set of support neighbors of $v$, and $N_p(v)$ is the set of all pendant neighbors of $v$.  Let $S$ be non-empty subset of the vertex set of $G$. Then the induced subgraph 
$\langle S \rangle$, is the graph with vertex set $S$ and  two vertices $u$ and $v$ are adjacent in $\langle S \rangle$ if and only if $u$ and $v$ are adjacent in $G$. If $T$ is tree of odd diameter, then $T$ has two central vertices $v_{c_1}$ and $v_{c_2}$ and if $T$ is tree of even diameter then $T$ has unique central vertex $v_c$. Note that none of the central vertices of $T$ is pendant when $diam(T)$ is $4$ or $5$.  

 Here, we recall certain well-known product graphs. For two graphs $G$ and $H$, the direct product, $G\times H$, is defined with vertex set $\{(g,h):g \in G, h\in H\}$ and edge set $\{((g,h)(g',h')):gg'\in E(G) \textnormal{~and~} hh'\in E(H)\}$.
 
 Frucht and Harary \cite{MR281659}
introduced the corona product of graphs. The corona product $G\odot H$, obtained by taking one copy of $G$ and $m$ copies of $H$ and joining each vertex from the $i$-th copy of $H$ with the $i$-{th} vertex of $G$ by an edge. 

For a graph $G$ with order $n$ and for given graphs $H_1, H_2, \ldots, H_n$, the generalised corona $G \tilde{\circ}\bigl(\displaystyle \mathop{\Lambda}\limits^n_{i=1} H_i\bigr)$, is the graph obtained by taking one copy of graphs $G, H_1, \ldots, H_n$ and joining the $i$-th vertex of $G$ to every vertex of $H_i$. 

The composition, $G[H]$, is defined with vertex set $\{(g,h):g \in G, h\in H\}$ and edge set $\{((g,h)(g',h')): g \sim g' \in G \textnormal{~or~} g=g' \textnormal{~and~} h\sim h' \in H\}$. It is also called the graph lexicographic product.
% \begin{definition}{\label{def5}}
% Harary graph $H_{k,n}$  is a $k$-connected graph on $n$ vertices has degree at least $k$
% and $\lceil{\frac{kn}{2}\rceil }$ edges. When $k$ is even the graph $H_{k,n+}$ is drawn by joining each vertex to its adjacent $\frac{k}{2}$ neighbors in each direction. When $k$ odd and $n$ even, take $k = 2r + 1$. Then $H_{2r+1,n}$ is constructed by first drawing $H_{2r,n}$ and then adding edges joining
% vertex $i$ to $i+ \frac{ n}{2}$ for $1 \le  i \le \frac{n}{2}$. Further when both $k$ and $n$ odd, Let $k = 2r + 1$. then $H_{2r+1,n}$ is constructed by first drawing $H_{2r,n}$ and then adding edges joining
% vertex $0$ to vertices $\frac{n-1}{2}$ and $\frac{n+1}{2}$
% and vertex $i$ to vertex $i + \frac{n+1}{2}$ for $1 \le  i < n-1$
% \end{definition}

If $G$ be a connected graph, then $V(G)$ can be partitioned in to pairwise disjoint sets $\{\Omega_{\mathfrak{p}}, \Omega_{\mathfrak{s}},  \Omega_{w}, \Omega_{\mathfrak{n}}\}$, where $\Omega_{\mathfrak{p}}$ is the set of all pendant vertices, $\Omega_{\mathfrak{s}}$ is the set of all strong support vertices, $\Omega_{w}$  is  the set of all  weak support  vertices, and $\Omega_{\mathfrak{n}}$ is the set of all neutral vertices. 
%-----------------
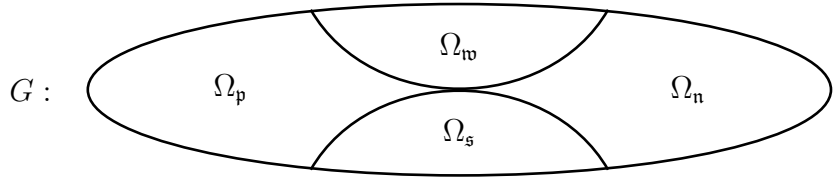
\begin{figure}[htbp]
	\centering

	\begin{tikzpicture}[scale=0.38]
 
  % \draw[rounded corners=30pt,rotate=180, line width=1pt] (0,0) ellipse (13cm and 3cm);
  % \draw[rounded corners=30pt,rotate=180, line width=1pt] (0,-3) circle (3cm and 3cm);
    % \draw[rounded corners=30pt,rotate=180, line width=1pt] (0,0) circle (5cm and 3cm);
  \draw[rounded corners=30pt,rotate=180, line width=1pt] (0,0) ellipse (13cm and 3cm);
   % \draw[line width=1pt] (-5,0)--(5,0);
  % \filldraw[fill=white,rounded corners=20pt,rotate=180, line width=1pt] (0,1.5) circle (5cm and 1.5cm); 
  %  \draw[fill=white,rounded corners=20pt,rotate=180, line width=1pt] (0,-1.5) circle (5cm and 1.5cm); 
 %  \draw[ rounded corners=30pt,rotate=180, line width=0.5pt] (0,0) ellipse (9cm and 2cm);
 % \draw[ rounded corners=30pt,rotate=180, line width=0.5pt] (0,0) ellipse (5cm and 1cm);
  \draw [line width=1pt] (-5.2,2.8) to [out=-60,in=-120,looseness=1.05] (5.2,2.8);
   \draw [line width=1pt] (-5.2,-2.8) to [out=60,in=120,looseness=1.05] (5.2,-2.8);
          \node at (-15,0){$G:$};	
  \node at (8,0){$\Omega_{\mathfrak{n}}$};	
   \node at (-8,0){$\Omega_{\mathfrak{p}}$};
      \node at (0,1.5){$\Omega_{\mathfrak{w}}$};
     \node at (0,-1.5){$\Omega_{\mathfrak{s}}$};	
%    \draw[line width=0.5pt] (0,1)--(0,2);
%    \draw [line width=0.5pt] (0,-1)--(0,-2);
\end{tikzpicture}
	
	\end{figure}
% The distance $d(u,v)$ between two vertices $u$ and $v$ is the length of a shortest $u-v$ path in $G$. The diameter of $G$ is defined by $diam(G)=\max\{d(u,v):u,v\in V\}$. The eccentricity of a vertex $v$ is the distance from $v$ to the vertex which is farthest from $v$. The center $C(G)$ of a graph $G$, is the set of vertices with minimum eccentricity. For a tree $T$, 
% $$C(T)=\begin{cases}\{v_c\} & \textnormal{if}~diam(T)~\textnormal{is~even}\\
% \{v_{c_1},v_{c_2}\} & \textnormal{if}~ diam(T) ~\textnormal{is~odd.} 
% \end{cases}$$
% Also, if $|C(T)|=2$, then the two central vertices $v_{c_1}$ and $v_{c_2}$ are  adjacent in $T$.  A bi-star $B_{r,s}$ is a tree with the vertex set $\{u,v, u_i,v_j: 1\le i\le r, 1\le j\le s\}$ and the edge set $\{uv, uu_i,vv_j: 1\le i\le r, 1\le j\le s\}$. Note that  $B_{r,s}$ has two strong support vertices $u$ and $v$, when $r\ge 2, s\ge 2$. 

\par Throughout this paper, $\Gamma$ denotes an additive Abelian group with the identity zero and $|\Gamma|$ denotes the order of  $\Gamma$. The order of an element $g$ of $\Gamma$ is denoted by $o(g)$. An element $g$ of $\Gamma$ is an involution if $g=-g$, where $-g$ is additive inverse of $g$ in $\Gamma$. 

An element $g$ of an Abelian group is a torsion element if $o(g)$ is finite.
If all of the elements of an Abelian group are torsion elements, the group is a torsion group. An Abelian group is torsion-free if it has no torsion elements other than zero. Note that in a torsion-free Abelian group, the identity element is the only element having finite order. Also note that, for any Abelian group $\Gamma$, the  set of all torsion elements $\Gamma_T$ forms a subgroup of $\Gamma$.

Let $\Gamma$ be a group with the identity element 0 and let $k \geq  1$ be a positive integer. We say that $\Gamma$ has exponent $n$ if $n$ is the smallest positive integer such that $\underbrace{g + \cdots + g}_{n~ \text{times}}=0$, for every $g\in \Gamma$. We use $e(\Gamma)$ to denote the exponent of $\Gamma$.  
If for a given group $\Gamma$,  there is no such
positive integer $k$, the convention is to set $e(\Gamma) =\infty$. Note that if $e(\Gamma) <\infty$, then the group $\Gamma$ is torsion and the converse is not. Consider the Prüfer $p$-group $\mathbb{Z}(p^\infty)=\cup_{t=1}^\infty \mathbb{Z}(p^t)$, where $p$ is a prime and $\mathbb{Z}(p^t)$ is the group of the $p^t$-th roots of unity. The group $\mathbb{Z}(p^\infty)$  is a torsion group with  $e(\mathbb{Z}(p^\infty))=\infty$. If $\Gamma$ is a group such that $e(\Gamma)$ is finite, then $e(\Gamma)$ is the least common multiple
of the orders of all the elements of $\Gamma$. In addition, if $\Gamma$ is a finite group, then the exponent $e(\Gamma)$ divides $|\Gamma|$.

%For a torsion group  $\Gamma$, the exponent of $\Gamma$, $e(\Gamma)$ is defined by \[\textcolor{red}{e(\Gamma)=\min\{k \in \mathbb{Z}^+ : kg= \underbrace{g + \cdots + g}_{k~ \text{times}}=0, \textnormal{~for all~} g \in \Gamma\setminus\{0\}\} }\] 
%is the smallest positive integer $k$ such that  $kg=\underbrace{g + \cdots + g}_{k~ \text{times}}=0$, for all $g\ne 0$. 
Moreover, for an Abelian group $\Gamma$, which is not torsion-free, we have $e(\Gamma)=e(\Gamma_T)$.  For group theoretic terminology, we refer to Herstein \cite{MR0356988}.
%\begin{proposition}{\label{prop1}}
%Every group of even order contains an involution. 
%\end{proposition}
\begin{definition}{\label{def1.2}}
The direct product of two groups  $(H,\circ)$ and $(K,\ast)$ is group $H\times K$ with  $(h_{1},k_{1})\circledast (h_{2}, k_{2})= (h_{1}\circ h_{2}, k_{1}\ast k_{2})$. Further, if $H$ and $K$ are Abelian groups, then $H\times K$ is an Abelian group so.
\end{definition}

A magic square is a $n\times n$ array in which the elements $1,2,\cdots,n^2$ appear exactly once and add up to the same sum in any row, column, main diagonal, or main back diagonal. Several authors proposed labelings that generalise the concept of magic squares. The reader can find an excellent treatment of various types of magic labeling in \cite{MR3013201}, \cite{MR1668059}. %Wallis and Marr's \cite{MR3013201} book. 
In 2001, Lee et al. \cite{MR1853016} introduced the concept of group-magic graphs. 
% \begin{definition}\textnormal{\cite{MR1853016}}
%     Let $\Gamma$ be an Abelian group. A graph $G = (V, E)$ is said  to be $\Gamma$-magic if there exists a labeling $\ell : E\rightarrow \Gamma$ such that the induced vertex
%  labeling $\ell^+ : V \rightarrow \Gamma$ defined by
% \[\ell^+(v) = \sum_{uv \in E} \ell (uv)\] is a constant map.
% \end{definition}
In 2020, Kamatchi et al.\cite{MR4145437}, introduced and studied group vertex magic graphs.
% In 2019, Kamatchi et al. \cite{MR4145437}, introduced and studied the concept of group vertex magic graphs. 
\begin{definition}\textnormal{\cite{MR4145437}}
If there is an element of $\mu$ of $\Gamma$ such that  $w(v)=\sum_{u\in N(v)}$ $\ell (u)=\mu$ for any vertex $v$ of $G$, then the mapping $\ell :V\rightarrow \Gamma$  is said to be an $\Gamma$-vertex magic labelling of $G$. A graph $G$ that allows such labelling is known as an $\Gamma$-vertex magic graph, and the corresponding $\mu$ is known as a magic constant. If $G$ has an $\Gamma$-vertex for each non-trivial Abelian group $\Gamma$, it is a group vertex magic graph. 
\end{definition} 
% \st{Note that if $G$ and $H$ are $\Gamma$-group vertex magic graphs, then $\mu(G)$ and $\mu(H)$ denote their respective magic constants.}
\newpage 
In 2022, Balamoorthy et al.\cite{MR4517978} discussed group vertex magicness of join
and tensor product of graphs and some basic results. 
\begin{theorem}{\label{1.5}}\textnormal{\cite{MR4517978}}
Let $\ell$  be an $\Gamma$ -vertex magic labeling of a graph $G$. Then
$\displaystyle{\sum_{v\in V(G)}\deg (v)\ell(v)= n\mu}$, where $n$ denotes the number of vertices of $G$ and $\mu$ is the magic constant.
\end{theorem}
 Recently, Sabeel et al.\cite{MR4515472} characterized the group vertex magicness of trees up to diameter $5$ for all finite Abelian groups.  In section {\ref{sec.3}}, we will discuss the group vertex magic of trees with diameter $4$ and $5$ for infinite Abelian groups having finitely many torsion elements. The following theorems are useful for proving our main results. 
 
\begin{lemma}\label{lemma1.6} \textnormal{\cite{MR4515472}} Let $\Gamma$ be an  Abelian group with $|\Gamma|\ge 3$ and let $g\in \Gamma$.  Then for each $n\ge 2$,  there exist  $g_1,g_2,\dots, g_n$ in $\Gamma\setminus\{0\}$  such that  $g=  g_1+g_2+\cdots+g_n$.
\end{lemma} 
\begin{theorem}
 \label{diam4} \textnormal{\cite{MR4515472}}
 	Let $\Gamma$ be a finite Abelian group with $|\Gamma|\ge 3$ and $T$ be a diameter $4$ tree with  central vertex $v_c$. Then $T$ is  $\Gamma$-vertex magic if and only if $T$ satisfies one of the criteria below:\\
 		\textnormal{(i)} Any non-pendant vertex of $T$ is in $\Omega_{\mathfrak{s}}(T)$\\
 		\textnormal{(ii)}  $v_c\in \Omega_{\mathfrak{s}}(T)$, $\deg(v_c)\not\equiv 2 \pmod {e(\Gamma)}$, and all the other all other non-pendant vertices are in $\Omega_{\mathfrak{s}}(T)$.\\
		\textnormal{(iii)}  $v_c\in \Omega_{\mathfrak{n}}(T)$ and $\gcd(\deg(v_c)-1,|\Gamma|)\neq 1.$
 \end{theorem}

\begin{theorem}\label{internal}\textnormal{\cite{MR4515472}}
 Let $\Gamma$ be an Abelian group with $|\Gamma|\ge 3$. If $G$ is a graph in which every non-pendant vertex is in $\Omega_{\mathfrak{s}}(G)$, then $G$ is $\Gamma$-vertex magic.
\end{theorem}  
\begin{theorem}{\label{01}}\textnormal{\cite{MR4515472}}
 	Let $\Gamma$ be a finite Abelian group such that $|\Gamma|\ge 4$ and $T$ be a  diameter $5$ tree such that  $v_{c_{1}}\notin \Omega_{\mathfrak{n}}(T)$ and $v_{c_{2}}\in \Omega_{\mathfrak{n}}(T)$. Then $T$ is $\Gamma$-vertex magic if and only if $T$ satisfies one of the criteria below:\\
		\textnormal{(i)} $\gcd(\deg(v_{c_2})-1,|\Gamma|)\ne 1$, \\
		\textnormal{(ii)} $v_{c_1}$ is not adjacent to a weak support vertex.
  \end{theorem}
  \newpage
\section{Main results}
In this section, we prove some sufficient conditions of group vertex magic graphs. Further, we prove certain necessary and sufficient conditions of graphs with  $|\Omega_p|\geq 1$. If $G$ is $\Gamma$-vertex magic graph, then $\Sigma_{\ell(G)}$ denotes the label sum of all the vertices of $G$ under the given labeling $\ell$. Moreover, we consider only finite Abelian groups in this section.
\begin{theorem}
    Suppose $G$ is a  group vertex magic graph with labeling $\ell$ and magic constant $\mu$ such that $\Sigma_{\ell(G)}=\mu$. Then $G$ is not a complete graph.
\end{theorem}
\begin{proof}
    Let $G$ be a complete graph with vertices $u_1,u_2, \ldots, u_n$. Suppose $G$ is a group vertex magic graph with $\Sigma_{\ell(G)}=\mu$. Since $G$ is complete, $w(u_1)=\mu=\ell(u_2)+\ell(u_3)+\cdots+\ell(u_n)$. Thus, \begin{align*}
        \ell(u_1)+\ell(u_2)+\cdots+\ell(u_n) & =\mu,
    \end{align*}
    which implies $\ell(u_1)=0$, a contradiction.
\end{proof}
\begin{theorem}{\label{thm2.1}}
Suppose a graph $G$ has two vertices $u$ and $v$ with $|N(u)|-1=|N(v)|$. If $G$  has a unique vertex $x$ such that $x\in N(u)$ and $x\not\in N(v)$, then $G$ is not group vertex magic. 
\end{theorem}
\begin{proof}
Let $x$ be a unique vertex in $G$ such that $x \in N(u)$ and $ x \notin N(v)$. Then $N(u) \setminus N(v)=\{x\}$. Suppose $G$ is a group vertex magic with labeling $\ell$. Then $w(v)=w(u)$ for $u$ and $v$, that is,

 \begin{align*}
 \sum_{y \in N(v)}\ell(y) &  = \sum_{y \in N(u)}\ell(y)\\
    \sum_{y \in N(v)}\ell(y) & = \sum_{y \in N(u)\setminus \{x\}}\ell(y)+\ell(x) \\  & \Rightarrow   \sum_{y \in N(v)}\ell(y) - \sum_{y \in N(u)\setminus \{x\} }\ell(y) = \ell(x).
\end{align*}
Hence $\ell(x)=0$, a contradiction. Thus, $G$ is not a group vertex magic.
\end{proof}

\begin{theorem}{\label{thm2.2}}
If $G$ is any graph having two vertices $u$ and $v$
such that $\deg(u)=1,\deg(v)=2$ and $d(u,v)=2$, then $G$ is not group vertex magic.
\end{theorem}
\begin{proof}
\noindent Let $\Gamma$ be any Abelian group with   $|\Gamma|\ge 3$. Assume that $\ell$ be a $\Gamma$-vertex magic labeling of $G$ with $\mu=g$, where $g\in \Gamma $. Let $u$ be a pendant vertex in $N(N(v))$. Then $w(u)= \ell(x)$, where $x \in N(v)$, which implies $\ell(x)=g$, hence $g \neq 0.$ 
Since $d(u,v)=2$, we have $w(v)=\ell(x)+\ell(y)$ implies $g=g+\ell(y)$. Hence, $\ell(y)=0,$
a contradiction.
	
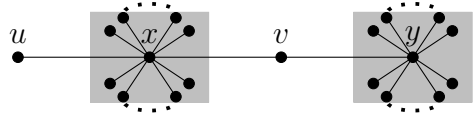
\begin{figure}[htbp]
 	\centering
 	\begin{tikzpicture}[scale=0.35]
          \filldraw[lightgray,line width=45pt,rounded corners=4pt] (2,-2)--(2,2)--cycle;
          \filldraw[lightgray,line width=45pt,rounded corners=4pt] (12,-2)--(12,2)--cycle;
        
 	\filldraw[fill = black, draw = black] (-3,0) circle (0.2 cm);
 	\filldraw[fill = black, draw = black] (2,0) circle (0.2 cm);
 	\draw (-3,0)--(2,0);
			\draw [line width=0.5mm, loosely dotted] (3,1.5) to [out=60,in=110,looseness=1] (1,1.5);
 \draw [line width=0.5mm, loosely dotted] (3,-1.5) to [out=-60,in=-110,looseness=1] (1,-1.5);
 
 	\filldraw[fill = black, draw = black] (3.5,-1) circle (0.2 cm);
 	\filldraw[fill = black, draw = black] (3.5,1) circle (0.2 cm);
 	\filldraw[fill = black, draw = black] (0.5,-1) circle (0.2 cm);
 	\filldraw[fill = black, draw = black] (0.5,1) circle (0.2 cm);
 	\draw (2,0)--(3.5,-1);
 	\draw (2,0)--(3.5,1);
  \draw (2,0)--(0.5,-1);
 	\draw (2,0)--(0.5,1);
 					\filldraw[fill = black, draw = black] (3,1.5) circle (0.2 cm);
      \filldraw[fill = black, draw = black] (1,1.5) circle (0.2 cm);
      \filldraw[fill = black, draw = black] (1,-1.5) circle (0.2 cm);
 				\draw (2,0)--(3,1.5);
     \draw (2,0)--(1,1.5);
     \draw (2,0)--(1,-1.5);
 	\filldraw[fill = black, draw = black] (3,-1.5) circle (0.2 cm);
 				\draw (2,0)--(3,-1.5);
\draw (12,0)--(2,0);
 \filldraw[fill = black, draw = black] (12,0) circle (0.2 cm);
  \filldraw[fill = black, draw = black] (7,0) circle (0.2 cm);

  \filldraw[fill = black, draw = black] (13.5,-1) circle (0.2 cm);
 	\filldraw[fill = black, draw = black] (13.5,1) circle (0.2 cm);
 	\filldraw[fill = black, draw = black] (10.5,-1) circle (0.2 cm);
 	\filldraw[fill = black, draw = black] (10.5,1) circle (0.2 cm);
  \filldraw[fill = black, draw = black] (11,1.5) circle (0.2 cm);
      \filldraw[fill = black, draw = black] (11,-1.5) circle (0.2 cm);
       \filldraw[fill = black, draw = black] (13,1.5) circle (0.2 cm);
      \filldraw[fill = black, draw = black] (13,-1.5) circle (0.2 cm);
       \draw (12,0)--(11,1.5);
   \draw (12,0)--(11,-1.5);
     \draw (12,0)--(13,1.5);
       \draw (12,0)--(13,-1.5);
         \draw (12,0)--(10.5,1);
           \draw (12,0)--(10.5,-1);
             \draw (12,0)--(13.5,1);
               \draw (12,0)--(13.5,-1);
 	\node at (-3,0.8) { $u$};
 	
 	\node at (7,0.8) { $v$};
 		\node at (12,0.8) { $y$};
 			\node at (2,0.8) { $x$};
    \draw [line width=0.5mm, loosely dotted] (13,1.5) to [out=60,in=110,looseness=1] (11,1.5);
 \draw [line width=0.5mm, loosely dotted] (13,-1.5) to [out=-60,in=-110,looseness=1] (11,-1.5);
 
 	\end{tikzpicture}
 	
 \end{figure}
\end{proof}

  %---------------
  \begin{corollary}{\label{lem2.5}} 
If $G$ contains a non-pendant vertex $u$ and  a pendant vertex $v$ such that $d(u,v)=2$, then $\deg(u)\geq 3$.  
  \end{corollary}
  \begin{proof}
Suppose $u$ is a vertex degree $2$, and $v$ is a pendant vertex in $G$ such that the distance between $u$ and $v$ is $2$. Then by Theorem {\ref{thm2.2}}, $G$ is not a group vertex magic graph.
 \end{proof}
  
  %--------------------------------
 
  \begin{theorem}
 Let $\Gamma$ be any Abelian group with $|\Gamma|\ge 3$ and $g\in \Gamma\setminus\{0\}$. Let $G$ be a graph in which every non-pendant vertex is weak support. Then $G$ is $\Gamma$-vertex magic with magic constant $g$ if and only if for any $x\in \Omega_w$,  $\deg(x)\not\equiv 2 \pmod {o(g)}$.

  \end{theorem}
  \begin{proof}
  Let $\Gamma$ be any Abelian group with $|\Gamma|\ge 3$. Assume that $G$ is a $\Gamma$-vertex magic graph with labeling $\ell$ and $\mu=g$, where $g\in \Gamma$.  Let $x$ and  $v$ be any two adjacent vertices in $G$ such that $x \in \Omega_{w}$  and $v \in \Omega_{\mathfrak{p}}$. Then $w(v)=\ell(x)=g$ and hence $g\neq 0$. Now,
 \[g=w(x)= \ell(v)+\sum_{x\in N(x)\setminus \{v\}} \ell(x)= \ell(v)+ (\deg(x)-1)g. \]
  Thus, $\ell(v)= (2-\deg(x))g$. Since $\ell(v) \neq 0$, we have $\deg(x) \not \equiv 2 \pmod{o(g)}$. 
 \\ For converse, let $g$ be any non-zero element of $\Gamma$. Consider a mapping $\ell: V(G) \rightarrow \Gamma$  such that \[\ell(z)=
  \begin{cases}
  g & \textnormal{~if~} z = u\\
  (2-\deg(x))g & \textnormal{~if~} z = v, \\
  \end{cases}
  \]
   \noindent where $x\sim v$ in $G$ such that $x \in \Omega_{w}$  and $v \in \Omega_{\mathfrak{p}}$. Since $o(g)$ does not divide $(2-\deg(x))$, $(2\deg(x))g\neq 0$. Then $w(x)=(\deg(x)-1)g+(2-\deg(x))g=g$, and $w(v)=g$. This completes the proof.
 \end{proof}
 Now, suppose $G$ is any graph with $n$ vertices in which all non-pendant vertices are in $\Omega_s$. Then the group vertex magicness of such $G$ is similar to find the group vertex magicness of generalized corona product of $G' \tilde{\circ} \bigl(\mathop{\Lambda}\limits^n_{i=1}H_i\bigr)$,  where $|V(G')|=n$ and each $H_i$ is the complement of complete graphs $K_m , m\ge 2$. 
  %------------------------------
\begin{theorem}  Let $\Gamma$ be any Abelian group with $|\Gamma|\ge 3$. Let $G$ be a graph with $n$ vertices, where $n>1$. Then the generalised corona product $G\tilde{\circ} \bigl(\mathop{\Lambda}\limits^n_{i=1}\overline{K}_{m_i}\bigr)$ is $\Gamma$-vertex magic if and only if $\deg(x)\geq3 $ for any $x\in V(G)$.
\end{theorem}
\begin{proof}
Let $G^\dagger= G\tilde{\circ} \bigl(\mathop{\Lambda}\limits^n_{i=1}\overline{K}_{m_i}\bigr)$. In $G^\dagger$, every non-pendant vertex is in $\Omega_\mathfrak{s}$. Therefore, every non-pendant vertices is adjacent with at least two pendant vertices. Assume that $G^\dagger$ is a group vertex magic graph with $\mu=g$. Suppose $u\in V( G)$ with $\deg(u)\le 2$. Then there exists a  vertex $v \in \overline{K}_{m_i}$ such that $d(u,v)=2$, which is a contradiction. Therefore, all vertices in $G$ have at least 3 degree. \\ 
For the converse, let the labeling $\ell$ from $V( G^\dagger)$ to $\Gamma$ such that label all the vertices of $G$ by any non-zero element $g$ of $\Gamma$ and label the vertices of $\overline{K}_{m_i}$ in such a way that their label sum is equal to $((m_i+1) - \deg(u_i))g$, by using Lemma {\ref{lemma1.6}}. Then,  weight of every pendant vertex is $g$, and
% \begin{align*}
% w_{\langle \Omega_{\mathfrak{p}} (G^\dagger)\rangle} (u) &= g, \textnormal{~~for all pendant vertices of } ~G^\dagger
% \end{align*} 
for any $i \in \{1,2,\cdots, n\}$, the weight of $u_i$
 \begin{align*}
 w_{G^\dagger} (u_i) &= (\deg(u_i)-m_i)g+((m_i+1)\deg(u_i))g=g.
 \end{align*}
This completes the proof.
\end{proof}

  %---------------------------------
  \begin{example}
      \noindent Consider the graph $G$ with vertices $u_1,u_2,\cdots u_8$ and $E(G)=\{u_iu_{i+1}, u_8u_1, u_1u_5, u_3u_7: 1\leq i\leq 7\}$. Construct the new graph \\ $G^\dagger=G\tilde{\circ}(\overline{K}_3,\overline{K}_4,\overline{K}_4,\overline{K}_4,\overline{K}_5,\overline{K}_6,\overline{K}_4,\overline{K}_5,\overline{K}_5)$. Suppose $\ell$ a mapping from $V(G^\dagger)$ to any Abelian group $\Gamma$ such that $\ell(u_i)=g$ and label the remaining vertices in such a way that $\Sigma_{\ell(\overline{K}_{m_i})}=((m_i+1)-\deg(u_i))g$. Clearly, the weight of any vertex of $G^\dagger$ is $g$.
  \end{example}
  \begin{figure}[htbp]
	\centering
	\begin{tikzpicture}[scale=0.25]
   %  \filldraw[ draw = black] (0,8) ellipse (1 cm);
  \filldraw[draw=lightgray, fill=lightgray,rounded corners=30pt,rotate=180, line width=1.5pt] (-0.25,8) ellipse (3.8cm and 1.2cm);
   \filldraw[draw=lightgray, fill=lightgray,rounded corners=30pt,rotate=180, line width=1.5pt] (0,-8) ellipse (3cm and 1.2cm);
    \filldraw[draw=lightgray, fill=lightgray,rounded corners=30pt,rotate=90, line width=1.5pt] (0,8) ellipse (3.2cm and 1.2cm);
\filldraw[draw=lightgray, fill=lightgray,rounded corners=30pt,rotate=90, line width=1.5pt] (0.2,-8) ellipse (3cm and 1.2cm);
  \filldraw[draw=lightgray, fill=lightgray,rounded corners=80pt,rotate=40, line width=1.5pt] (1.2,-12.8) ellipse (4.2cm and 1.2cm);
   \filldraw[draw=lightgray, fill=lightgray,rounded corners=80pt,rotate=40, line width=1.5pt] (-1.2,12.8) ellipse (4.2cm and 1.2cm);
  \filldraw[draw=lightgray, fill=lightgray,rounded corners=80pt,rotate=-40, line width=1.5pt] (-1.3,-12.6) ellipse (3.8cm and 1.2cm);
   \filldraw[draw=lightgray, fill=lightgray,rounded corners=80pt,rotate=-40, line width=1.5pt] (1.3,12.6) ellipse (3.8cm and 1.2cm);
   \draw[ draw = black] (0,0) circle (4 cm);
  \filldraw[ draw = black] (0,4) circle (0.2 cm);
   \filldraw[ draw = black] (0,-4) circle (0.2 cm);
    \filldraw[ draw = black] (4,0) circle (0.2 cm);
     \filldraw[ draw = black] (-4,0) circle (0.2 cm);
   \draw (4,0)--(-4,0);
   \draw (0,4)--(0,-4);
       \filldraw[ draw = black] (1,8) circle (0.2 cm); \filldraw[ draw = black] (0,8) circle (0.2 cm);
       \filldraw[ draw = black] (-1,8) circle (0.2 cm);   
         \filldraw[ draw = black] (1,-8) circle (0.2 cm); \filldraw[ draw = black] (0,-8) circle (0.2 cm);
       \filldraw[ draw = black] (-1,-8) circle (0.2 cm);   
         \filldraw[ draw = black] (2,-8) circle (0.2 cm); \filldraw[ draw = black] (3,-8) circle (0.2 cm);
       \filldraw[ draw = black] (-2,-8) circle (0.2 cm);   
        
      \filldraw[ draw = black] (8,0) circle (0.2 cm); \filldraw[ draw = black] (8,1) circle (0.2 cm);
       \filldraw[ draw = black] (8,2) circle (0.2 cm);   \filldraw[ draw = black] (8,-1) circle (0.2 cm);
       \filldraw[ draw = black] (-8,-2) circle (0.2 cm);
       \filldraw[ draw = black] (-8,0) circle (0.2 cm); \filldraw[ draw = black] (-8,1) circle (0.2 cm);
       \filldraw[ draw = black] (-8,2) circle (0.2 cm);   \filldraw[ draw = black] (-8,-1) circle (0.2 cm);
       \filldraw[ draw = black] (-8,-2) circle (0.2 cm);
 \draw (0,4)--(0,8);
  \draw (0,4)--(1,8);
   \draw (0,4)--(-1,8);
   \draw (4,0)--(8,0);
   \draw (4,0)--(8,1);
   \draw (4,0)--(8,2);
   \draw (4,0)--(8,-1);
   \draw (-4,0)--(-8,1);
   \draw (-4,0)--(-8,2);
   \draw (-4,0)--(-8,-1);
    \draw (-4,0)--(-8,-2);
 \draw (-4,0)--(-8,0);
  \draw (0,-4)--(0,-8);
  \draw (0,-4)--(1,-8);
   \draw (0,-4)--(2,-8);
    \draw (0,-4)--(3,-8);
     \draw (0,-4)--(-2,-8);
   \draw (0,-4)--(-1,-8);
   \node at (0,-10.5){$\overline{K}_6$};
    \node at (0,10.5){$\overline{K}_3$};
   \node at (11,0){$\overline{K}_4$};
    \node at (-11,0){$\overline{K}_5$};
     \node at (-12,11){$\overline{K}_5$};
     \node at (12,11){$\overline{K}_4$};
     \node at (-12,-11){$\overline{K}_4$};
     \node at (12,-11){$\overline{K}_5$};
 \node at (1,4.5){$u_1$};
  \node at (4.2,2.5){$u_2$};
 \node at (4.5,-1){$u_3$};
 \node at  (2.4,-4){$u_4$};
     \node at (-4.2,-2.5){$u_6$};
     \node at (-1.5,-4.5){$u_5$};
    \node at (-4.5,1.2){$u_7$}; 
    \node at (-2.5,3.9){$u_8$};
      \filldraw[ draw = black] (2.8,2.8) circle (0.2 cm);  
      
       \filldraw[ draw = black] (-2.8,-2.8) circle (0.2 cm);  
     
        \filldraw[ draw = black] (2.8,-2.8) circle (0.2 cm); 
       
         \filldraw[ draw = black] (-2.8,2.8) circle (0.2 cm); 
        
          \filldraw[ draw = black] (10,8) circle (0.2 cm);  
           \filldraw[ draw = black] (9,9) circle (0.2 cm);  
            \filldraw[ draw = black] (11,7) circle (0.2 cm);  
           \filldraw[ draw = black] (8,10) circle (0.2 cm);  
           \draw (2.8,2.8)--(10,8);
            \draw (2.8,2.8)--(8,10);
             \draw (2.8,2.8)--(11,7);
              \draw (2.8,2.8)--(9,9);

               \filldraw[ draw = black] (-10,-8) circle (0.2 cm);  
           \filldraw[ draw = black] (-9,-9) circle (0.2 cm);  
            \filldraw[ draw = black] (-11,-7) circle (0.2 cm);  
           \filldraw[ draw = black] (-8,-10) circle (0.2 cm);  
           \draw (-2.8,-2.8)--(-10,-8);
            \draw (-2.8,-2.8)--(-8,-10);
             \draw (-2.8,-2.8)--(-11,-7);
              \draw (-2.8,-2.8)--(-9,-9);

               \filldraw[ draw = black] (-7,11) circle (0.2 cm);  
          \filldraw[ draw = black] (-10,8) circle (0.2 cm);  
           \filldraw[ draw = black] (-9,9) circle (0.2 cm);  
            \filldraw[ draw = black] (-11,7) circle (0.2 cm);  
           \filldraw[ draw = black] (-8,10) circle (0.2 cm);  
           \draw (-2.8,2.8)--(-10,8);
            \draw (-2.8,2.8)--(-8,10);
             \draw (-2.8,2.8)--(-11,7);
              \draw (-2.8,2.8)--(-9,9);
\draw (-2.8,2.8)--(-7,11);
               \filldraw[ draw = black] (7,-11) circle (0.2 cm);  
          \filldraw[ draw = black] (10,-8) circle (0.2 cm);  
           \filldraw[ draw = black] (9,-9) circle (0.2 cm);  
            \filldraw[ draw = black] (11,-7) circle (0.2 cm);  
           \filldraw[ draw = black] (8,-10) circle (0.2 cm);  
           \draw (2.8,-2.8)--(10,-8);
            \draw (2.8,-2.8)--(8,-10);
             \draw (2.8,-2.8)--(11,-7);
              \draw (2.8,-2.8)--(9,-9);
\draw (2.8,-2.8)--(7,-11); 
\end{tikzpicture}
\caption{$G\tilde{\circ}(\overline{K}_3,\overline{K}_4,\overline{K}_4,\overline{K}_4,\overline{K}_5,\overline{K}_6,\overline{K}_4,\overline{K}_5,\overline{K}_5) $}
 \end{figure}
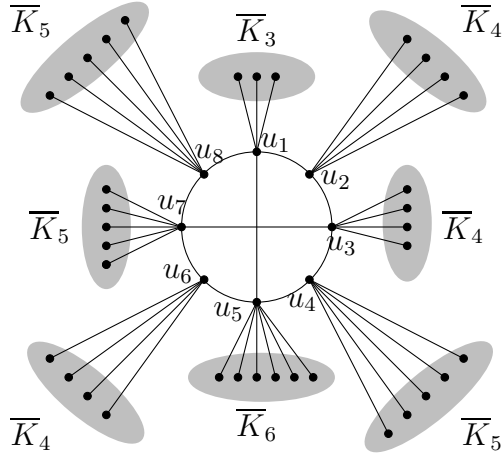

 %--------------------
 
\newpage In the following theorems, we assume that $x_i$'s are elements of $\Omega_{\mathfrak{n}}$, $y_i$'s are elements of $ \Omega_{\mathfrak{s}}\cup\Omega_{\mathfrak{w}}$, and $z_i$'s are elements of $\Omega_{\mathfrak{p}}$ in a graph $G$.
\begin{theorem}{\label{Thm2.8}}
   Let $\Gamma$ be any Abelian group with $|\Gamma|\ge 3$. Let $G$ be a graph with at least one pendant and all non pendant vertices are either in $\Omega_{\mathfrak{s}}$ or in $\Omega_{\mathfrak{n}}$.  If there exists an non-zero  element $g \in \Gamma$ such that $\deg(x_i)\equiv 1 \pmod{o(g)}$, for all $x_i \in \Omega_{\mathfrak{n}}$, then $G$ is $\Gamma$-vertex magic graph. 
\end{theorem}
\begin{proof} Assume that $g$ is a non zero element in $g$ such that $\deg(x_i)\equiv 1 \pmod{o(g)}$, for all $x_i \in \Omega_{\mathfrak{n}}$.
   Consider a labeling $\ell:V(G)\rightarrow \Gamma$, such that label of all non pendant vertices in $G$ is $g$ and using Lemma {\ref{lemma1.6}}, label the pendant vertices adjacent to each $y_i's$ in such a way that their label sum is $(1-\deg_s(y_i)+\deg_n(y_i))g$. Then,
    % {Case.1.} If $x_i\sim x_j$ for all $i$ and $j$, and each neutral vertex is adjacent to atleast one strong support vertex. Then, 
    \begin{align*}
        w_G(z_i)&=g, \text{~for all~} i,\\
        w_G(y_i)&=\deg_s(y_i)g+\deg_n(y_i)g+(1-\deg_s(y_i)+\deg_n(y_i))g =g, \text{~for all~} i,\\
            w_G(x_i)&=\deg_s(x_i)g+\deg_n(x_i)g=\deg(x_i)g=g, \text{~for all~} i, 
\end{align*}
  because any $\deg(x_i)$ is congruent to 1 modulo $o(g)$. 
   % {Case.2.} Suppose $x_i\nsim x_j$ for all $i$ and $j$, and each neutral vertex is adjacent to atleast one strong support vertex. Then, 
     % \begin{align*}
        % w_G(z_i)&=g, \text{~for all~} i,\\
%         w_G(y_i)&=\deg_s(u_i)g+\deg_n(u_i)g+(1-\deg_s(u_i)+\deg_n(u_i))g =g, \text{~for all~} i,\\
% w_G(x_i)&=\deg(x_i)g=g, \text{~for all~} i, \text{~since~}\deg(x_i)\equiv 1 \pmod{o(g)}. 
%     \end{align*}
% Therefore, $G$ is $\Gamma$-vertex magic graph.  \\
%    {Case.2.} Suppose  $x_i\nsim x_j$ for some $i$ and $j$, and each neutral vertex is adjacent to atleast one strong support vertex. 

% And the remaining cases two two cases $x_i\nsim x_j$ for all $i$ and $j$, and some neutral vertex is adjacent to at least one strong support vertex,  $x_i\nsim x_j$ for some $i$ and $j$, and \color{red}{some neutral vertex is adjacent to at least one strong support vertex, are we can analogously proved by comparing case.1 and case .2.} 
\noindent Therefore $G$ is $\Gamma$-vertex magic graph.
\end{proof}
  
  %---------------
  \begin{theorem}{\label{thm2.5}} 
Let $\Gamma$ be any Abelian group with $|\Gamma|\geq 3$ and let $g \in \Gamma$.
Let $G$ be a graph with atleast one pendant vertex and exactly one neutral vertex $x$. Then $G$ is  $\Gamma$-vertex magic with magic constant $g$ if and only if $\deg (x)\equiv1 \pmod{o(g)}$ and $\deg (x)\not\equiv 2 \pmod{o(g)}$. 
  \end{theorem}
  \begin{proof}
Assume that $G$ is a $\Gamma$- vertex magic graph with magic constant $g$. 
% Let $y_i \in \Omega_{\mathfrak{s}}\cup\Omega_{w}, v_i \in \Omega_{\mathfrak{p}}$, where $i\in\{1,2,\cdots ,|\Omega_{\mathfrak{s}}|+|\Omega_{w}|\}$. 
Then for all $y_i \in \Omega_{\mathfrak{s}}\cup\Omega_{w}$, $\ell(y_i)=g $, and hence $g\neq 0$. 
For the neutral vertex $x$,  $w(x)=\deg(x)g =g$. Thus $(\deg(x)-1)g =0$ implies $\deg (x)\equiv1 \pmod{o(g)}$.
In addition, if $\deg(x)=2$, then $g=0$, which is not possible and hence $\deg (x)\not\equiv 2 \pmod{o(g)}$.\\
%  Let $u_1$ is adjacent to a normal vertex $x$ and pendant vertex $v_1$, then $$w(u_1)= (\deg (u_1)-2)g+ \ell(v_1)+\ell(x)$$
% Suppose $\deg (u_1)=2$, then $g= \ell(v_1)+\ell(x) $
    For the converse, suppose  $y_i\in \Omega_{\mathfrak{s}}$. Define $\ell : V(G)\rightarrow \Gamma$ such that label of all non pendant vertices as $g$ and by Lemma {\ref{lemma1.6}}, label the pendant vertices  adjacent to each $y_i's$ in such a way that their label sum is $((\deg_p(y_i)+1)-\deg(y_i))g$. Therefore weight must be $g$. Then, 
     \begin{align*}
        w_G(z_i)&=g, \text{~for all~} i,\\
       % w_G(y_i)&=\deg_s(y_i)g+\deg_n(y_i)g+(1-\deg_s(y_i)+\deg_n(y_i))g =g, \text{~for all~} i,\\
            w_G(x)&=\deg(x)g=g. \textnormal{~Therefore~} \deg (x)\equiv1 \pmod{o(g)}\\
         w_G(y_i)&=\deg_s(y_i)g+\deg_n(y_i)g+ (2-\deg(y_i))g\\
      &= (\deg(y_i)-\deg_p(y_i))g + ((\deg_p(y_i)+1)-\deg(y_i))g=g, \text{~for all~} i.   
\end{align*}

    % For any $u_i \in \Omega_{w}$, define  $\ell(v_{i})=(2-\deg(u_i))g$. \\
\noindent  Now for any $y_i\in \Omega_{\mathfrak{w}}$, define the label for pendant vertices adjacent to each $y_i$ as $g', g'\neq g $. Then, 
     \begin{align*}
      w_G(y_i)&=\deg_s(y_i)g+\deg_n(y_i)g+ g'\\
      &g= (\deg(y_i)-1)g +g'
 \end{align*}
 Therefore, take $g'=(2-\deg(y_i))g$. 
 Hence $G$ is a $\Gamma$-vertex magic graph.
  \end{proof}

  % From Theorem {\ref{thm2.5}}, and Lemma {\ref{lemma1.6}}, the following observations are obvious. 

  % \begin{observation}
  %      Let $G$ be a graph with exactly one normal vertex and all others are either strong support vertex or pendant vertex, then $G$ is  group vertex magic with $\mu, \mu\in \Gamma$ if and only if $x$ is a normal vertex such that $\deg(x)>2$ and $\deg (x)\equiv1 \pmod{o(\mu)}$. 
  % \end{observation}
  % \begin{proof}
  % The proof is similar to Theorem {\ref{thm2.5}}. 
  % \end{proof}

  \begin{theorem}{\label{thm2.8}}
   Let $\Gamma$ be any Abelian group with $|\Gamma|\geq 3$ and $g \in \Gamma$. Let $G$ be a graph in which $| \Omega_{\mathfrak{n}}(G)|=2$. Then $G$ is $\Gamma$-vertex magic with the magic constant $g$ if and only if for any non pendant vertex $v\in G$,  $\deg (v)\equiv1 \pmod{o(g)}$ and $\deg(v)\not\equiv 2\pmod{o(g)}$. 
     \end{theorem}
 \begin{proof}
         Assume that $G$ is   $\Gamma$-vertex magic  with labeling  $\ell$ and magic constant $g$. Then the label of any vertex in $(\Omega_\mathfrak{s}\cup\Omega_\mathfrak{w})(G)$ must be $g$, and hence $g\ne 0$.\\ Let $x_1,x_2\in \Omega_{\mathfrak{n}}$.\\
           {{Case 1.}} Suppose $ x_1\sim x_2$. Then, 
        \begin{align*}
            w_G(x_1) & =(\deg(x_1)-1)) g+\ell(x_2)\\
            \ell(x_2) & = (2-\deg(x_1))g.
        \end{align*}
     Similarly, $\ell(x_1)  = (2-\deg(x_2))g $. Suppose $\deg(x_1)=2$, then $\ell(x_1)=0$, which is a contradiction. Therefore $\deg(x_i)\not\equiv 2\pmod{o(g)},~ \textnormal{for }i=1,2$. 
     For the vertices in $\Omega_{\mathfrak{s}}\cup \Omega_{w}$, we have the following four subcases. 
     
      \noindent {{Case 1.1.}} Let $y_i\in \Omega_{w}$. Suppose $y_i\sim x_1$, $y_i\nsim x_2$, and $z_i\sim y_i$. Then, 
      \begin{align*}
          w_G(y_i) &= (\deg(y_i)-2)g+\ell(x_1)+\ell(z_i)\\
          \ell(z_i) &= (\deg(x_2)-\deg(y_i)-1)g
      \end{align*}
      Similarly, $y_i\sim x_2$, $y_i\nsim x_1$, and $z_i\sim y_i$. Then $\ell(z_i)=(\deg(x_1)-\deg(y_i)-1)g$. \\
       {{Case 1.2.}} Let $y_i\in \Omega_{\mathfrak{s}}$. Suppose $z_i's$ are  the pendant vertices adjacent to $y_i$, and  $y_i \sim x_1$. Then, from Lemma {\ref{lemma1.6}},  label the pendant vertices in such a way that their label sum is  $(\deg(x_2)-\deg(y_i)-1)g$.\\
       Similarly, $y_i \sim x_2$. Using the Lemma {\ref{lemma1.6}}, label the pendant vertices so that their label sum is $(\deg(x_1)-\deg(y_i)-1)g$.\\
{{Case 1.3.}} Let $y_i \in \Omega_{w}$. Suppose $y_i\sim z_i$, $y_i\not \sim x_1$, and $y_i\nsim x_2$. Then,
\begin{align*}
          w_G(y_i) &= (\deg(y_i)-1)g+\ell(z_i),\\
          \ell(z_i) &= (2-\deg(y_i))g.
      \end{align*}
 Similarly, if $y_i$ is in $\Omega_{\mathfrak{s}}$, using lemma {\ref{lemma1.6}}, label the pendant vertices adjacent $y_i$ in such a way that their label sum is $(2-\deg(y_i))g$.\\
{{Case 1.4.}} Let $y_i \in \Omega_{w}$. Suppose $y_i\sim z_i$, $y_i \sim x_1$, and $y_i\sim x_2$. Then,
\begin{align*}
          w_G(y_i) &= (\deg(y_i)-3)g+\ell(x_1)+\ell(x_2)+\ell(z_i),\\
          g &= (\deg(y_i)-3)g+(2-\deg(x_2))g+(2-\deg(x_1))g+\ell(z_i),\\
          \ell(z_i) &= (\deg(x_1)+\deg(x_2)-\deg(y_i)-1)g.
      \end{align*}
 Similarly, if $y_i$ is in $\Omega_{\mathfrak{s}}$, using Lemma {\ref{lemma1.6}}, label the pendant vertices adjacent $y_i$ in such a way that their label sum is $(\deg(x_1)+\deg(x_2)-\deg(y_i)-1)g$.\\
  {{Case 2.}} Suppose $ x_1\not \sim x_2$.  Then, $w_G(x_1)=\deg(x_1)g=g$ implies $\deg(x_1)\equiv 1 \pmod{o(g)}$. Similarly, $\deg(x_2)\equiv 1 \pmod{o(g)}$.  Let $g',g''\in \Gamma$ such that $\ell(x_1)=g'$ and $\ell(x_2)=g''$.\\
 {{Case 2.1.}} Let $y_i \in \Omega_{\mathfrak{w}}$. Suppose $y_i \sim x_1, y_i \nsim x_2$, and $y_i\sim z_i$. Then, \begin{align*}
     w_G(y_i) & = g=(\deg(y_i)-2)g+\ell(x_1)+\ell(z_i)\\
     \ell(z_i) &= (3-\deg(y_i))g-g'
 \end{align*} Analogously, $\ell(y_i) = (3-\deg(y_i))g-g''$, whenever $y_i\sim x_2$ and  $y_i\not \sim x_1$.\\
 Now suppose $y_i\in \Omega_{\mathfrak{s}}$, using Lemma {\ref{lemma1.6}}, label the pendant vertices.\\
  {{Case 2.2.}} Let $y_i \in \Omega_{\mathfrak{w}}$. Suppose $y_i \nsim x_1, y_i \nsim x_2$, and $y_i\sim z_i$. Then,
 \begin{align*}
     w_G(y_i) & = g=(\deg(y_i)-1)g+\ell(z_i)\\
     \ell(z_i) &= (2-\deg(y_i))g
 \end{align*} 
 If $y_i \in \Omega_{\mathfrak{s}}$, then by Lemma {\ref{lemma1.6}}, label the pendant vertices adjacent to $y_i$ such that their label sum is $(2-\deg(y_i))g$.\\
 {{Case 2.3.}} Let $y_i \in \Omega_{\mathfrak{w}}$. Suppose $y_i \sim x_1, y_i \sim x_2$, and $y_i\sim z_i$. Then,
\begin{align*}
     w_G(y_i) & = g=(\deg(y_i)-3)g+\ell(x_1)+\ell(x_2)+\ell(z_i)\\
     \ell(z_i) &= (2-\deg(y_i))g-g'-g''
 \end{align*} 
 Similarly, for any $y_i \in \Omega_{\mathfrak{s}}$, use Lemma {\ref{lemma1.6}}, to label its pendant vertices. \\
 Hence it completes the proof.
     \end{proof}
   %--------------------------------
   \begin{observation}
    Let $G$ be a graph with at least one pendant vertex.  Let $G$ be a $\Gamma$-vertex magic graph with labeling $\ell$ and magic constant $g$. Then,
       \begin{enumerate}
           \item label of all  vertices in $\Omega_{\mathfrak{s}}\cup \Omega_{w}$, which are always $g$. Hence, $g$ is not the identity element of any Abelian group.
           \item Let  $x_1, x_2 \in \Omega_{\mathfrak{n}}$ such that  $x_1\sim x_2 $, then $\ell(x_1)-\ell(x_2)=(\deg(x_1)-\deg(x_2))(-g)$. 
           \item Let $x_1, x_2 \in \Omega_{\mathfrak{n}}$ such that  $x_1\nsim x_2 $, then $(\deg x_1, \deg x_2) \equiv1 \pmod{o(g)}$  and $(\deg x_1, \deg x_2) \not\equiv2 \pmod{o(g)}$
           
           % If  $y_i\nsim y_j \in \Omega_{\mathfrak{n}}$ for any $i \neq j$, then $\deg(y_i)\not\equiv 2\pmod{o(g)}$, and $\deg (y_i)\equiv1 \pmod{o(g)}$ and $\deg(y_j)\not\equiv 2\pmod{o(g)}$, and $\deg (y_j)\equiv1 \pmod{o(g)}$. 
           \item If $y_1, y_2  \in \Omega_{\mathfrak{s}}\cup \Omega_{w}$, then $\deg(y_1)\ge 3$ and $\deg(y_2)\ge 3$.
       \end{enumerate}
   \end{observation}

\begin{proof}
Assume that $G$ be a $\Gamma$- vertex magic graph with magic constant $\mu=g$, then $\ell(y)=g, \textnormal{~~for all~~} y \in \Omega_{\mathfrak{s}}\cup \Omega_{w}$. Let $x_{i} \in \Omega_{\mathfrak{n}} , 1 \le i \le n.$\\ Suppose all $x_{i}$ are adjacent.  
Then,\small{ \begin{align*}
    w(x_i) & = (\deg(x_i)-(n-1))g +\sum_{v \in N_{\Omega_{\mathfrak{n}}}(x_i)}\ell(v)\\
    g & = (\deg(x_i)-(n-1))g +\sum_{v \in N_{\Omega_{\mathfrak{n}}}( x_i)}\ell(v)\\
   (n- \deg(x_i))g  & = \sum_{v \in \Omega_{\mathfrak{n}}(x_i)}\ell(v)\\
    \end{align*}}
    That is,
    \small{  \begin{equation}{\label{eq*}}
  \left.\begin{aligned}
     (n- \deg(x_1))g & = \ell(x_2)+\ell(_3)+\ell(x_4)+\cdots +\ell(x_n)\\
     (n- \deg(x_2))g& = \ell(x_1)+\ell(x_3)+\ell(x_4)+\cdots +\ell(x_n)\\
   %  (k- \deg(x_3))g& = \ell(x_1)+\ell(x_2)+\ell(x_4)+\cdots +\ell(x_n)\\
     \vdots~~~~~ & ~~~~~~\vdots ~~~~~~~~~~ \vdots ~~~~~~~~~~\vdots ~~~~~~~~~~~~~~\vdots\\
     (n- \deg(x_n))g& = \ell(x_1)+\ell(x_2)+\ell(x_3)+\cdots +\ell(x_{n-1})\\
 \end{aligned}\right.
 \end{equation}}
From equation {\ref{eq*}}, it is clear that, \[\ell(x_i)-\ell(x_j)=(\deg(x_j)-\deg(x_i))g= (\deg(x_i)-\deg(x_j))(-g).\]
Case 2. is satisfied. It is not hard to verify other cases. \end{proof}
%-----------------

\begin{theorem}
   Let $\Gamma$ be any Abelian group with $|\Gamma|\geq 3$ and let $g \in \Gamma\setminus\{0\}$.  Let $G$ be a graph such that  the induced subgraph $\langle\Omega_{\mathfrak{n}}\rangle$ contains more than one vertices of $G$. If $\langle\Omega_{\mathfrak{n}}\rangle$ is  $\Gamma$-vertex magic graph with $\mu=g$
    and  if for each vertex $x$ of $\Omega_{\mathfrak{n}}$, $\deg_s(x)\equiv 0\mod o(g)$, then $G$ is $\Gamma$-vertex magic graph with $\mu=g$.
\end{theorem}
\begin{proof}
 Given that $\langle\Omega_{\mathfrak{n}}\rangle$ is  $\Gamma$-vertex magic graph. Then there exists  $\ell$ from the vertex set of $\langle\Omega_{\mathfrak{n}}\rangle$ to $\Gamma$ such that $w_{\langle\Omega_{\mathfrak{n}}\rangle}(x)= g, \textnormal{~for any~} x \in \Omega_{\mathfrak{n}}$. \\  Consider  $\ell^\dagger : V(G)\rightarrow \Gamma$ such that \[\small{
\ell^\dagger(x)=
\begin{cases}
    g & \textnormal{if}~~ x \in \Omega_{\mathfrak{s}} \cup \Omega_{w}\\
    \ell(x) &  \textnormal{if}~~ x \in \Omega_{\mathfrak{n}}. \\
 %    &  \textnormal{If}~~ x \in \Omega_{\mathfrak{p}} \\
\end{cases}}
\]
% From the graph, it is clear that for  any vertex $y\in \Omega_{\mathfrak{s}} \cup \Omega_{w}$, $\deg(y)= k_1+ k_2 +k_3$, where $k_1$ is the number of vertices adjacent to $y$ in $\Omega_{\mathfrak{s}} \cup \Omega_{w}$, $k_2$ is the number of vertices adjacent to $y$ in $\Omega_{\mathfrak{n}}$, and $k_3$ is the number of vertices adjacent to $y$ in $\Omega_{\mathfrak{p}}$. \\
For the vertices in $\Omega_{\mathfrak{p}}$, there are two cases.\\
 Case (i): Let $z\in \Omega_{\mathfrak{p}}$  is adjacent to a vertex $y_i \in \Omega_{w}$. Then \[\small{\ell^\dagger(v)= (1-\deg_s(y_i))g -  \sum_{\scriptscriptstyle{v \in N_{\langle\Omega_{\mathfrak{n}}\rangle}(x)}}\ell(v)}.\]
 Case (ii): Let $t\ge 2$,  $z_i\in \Omega_{\mathfrak{p}}, i \in \{1,2,\cdots, t\}$ is adjacent to $y_i \in \Omega_{\mathfrak{s}}$. Label the vertices $v_i$ such that their label sum is $ \small{\big[(1-\deg_s(y_i))g -  \sum_{\scriptscriptstyle{v \in N_{\langle\Omega_{\mathfrak{n}}\rangle}(x)}}\ell(v)\big]}$ by using Lemma \ref{lemma1.6}. \\ Now, consider the weights of each vertices in $G$ as follows.\\ 
 \noindent  For any $z_i \in  \Omega_{\mathfrak{p}}$, the $w^\dagger_{\Omega_{\mathfrak{p}}}(z_i)=g$. Suppose $x \in \Omega_{\mathfrak{n}}$. Then the weight of $x$, \begin{align*}
w^\dagger_{\scriptscriptstyle{\Omega_{\mathfrak{n}}}}(x) & = \ell^\dagger_{\scriptscriptstyle{v \in N_{\Omega_{\mathfrak{n}}}(x)}}(v)+\ell^\dagger_{\scriptscriptstyle {v \in N_{(\Omega_{\mathfrak{s}} \cup \Omega_{w})}(x)}}(v)\\
 & = w_{\langle\Omega_{\mathfrak{n}}\rangle}(x) + (g +g + \cdots +g)\\
 &= w_{\langle\Omega_{\mathfrak{n}}\rangle}(x) + \deg_s(x)g\\
w^\dagger_{\scriptscriptstyle{\Omega_{\mathfrak{n}}}}(x )& = g ~~\textnormal{for all } x \in  \Omega_{\mathfrak{n}}.
\end{align*}
Suppose $y_i \in \Omega_{\mathfrak{s}}\cup \Omega_{w}$.
Then the weight of $y_i$, \small{\begin{align*}
w^\dagger_{\scriptscriptstyle{\Omega_{\mathfrak{s}} \cup\Omega_{w}}}(y_i)& =\ell^\dagger_{\scriptscriptstyle {v \in N_{(\Omega_{\mathfrak{s}} \cup \Omega_{w})}(y_i)}}(v)+ \ell^\dagger_{\scriptscriptstyle{v \in N_{\Omega_{\mathfrak{n}}}(y_i)}}(v)+ \ell^\dagger_{\scriptscriptstyle{v \in N_{\Omega_{\mathfrak{p}}}(y_i)}}(v) \\
 &= (g+g+\cdots + g)+ \sum_{\scriptscriptstyle{y \in N_{\langle\Omega_{\mathfrak{n}}\rangle}(y_i)}}\ell(v)+ \big[ (1-\deg_s(y_i))g -  \sum_{\scriptscriptstyle{v \in N_{\langle\Omega_{\mathfrak{n}}\rangle}(y_i)}}\ell(v)\big]\\
 w^\dagger_{\scriptscriptstyle{\Omega_{\mathfrak{s}} \cup\Omega_{w}}}(y_i)& =g  ~~\textnormal{for all } y_i\in  \Omega_{\mathfrak{s}}\cup \Omega_{w}.
\end{align*}}
From the above cases, it is clear that $G$ is $\Gamma$-vertex magic with $\mu =g$.
\end{proof}
% %--------------------------
\begin{theorem}
    Let $G$ be a graph such that all vertices in $\Omega_n$ are non adjacent. Let $g\in \Gamma$, such that the number of adjacent vertices of each vertex of $\Omega_{\mathfrak{s}} \cup\Omega_{w}$ in $\Omega_{\mathfrak{n}}$ is an
integral multiple of $o(g)$. Then $G$ is a $\Gamma$- vertex magic graph.
     \end{theorem}
\begin{proof}
    Given that $G$ is a graph such that all vertices in $\Omega_n$ are non adjacent. Therefore, from Theorem \ref{thm2.8}, $\deg(x)\not\equiv 2\pmod{o(g)}$ and $\deg (x)\equiv1 \pmod{o(g)}$ whenever $x\in  \Omega_{\mathfrak{n}}$. \\
    Consider the labeling $\ell$ from $V(G)$ to $\Gamma$ such that label of all non pendant vertices $g$ and for pendant vertices, we have the following cases.
    % \[\ell(v)=\begin{cases}
    % g & \textnormal{if~} v \in \Omega_{\mathfrak{s}} \cup\Omega_{w},  \textnormal{~or~} x \in \Omega_{\mathfrak{n}}\\
    
    % (1-\deg_s(x))g & \textnormal{If~} x \in \Omega_{\mathfrak{w}}.
    % \end{cases}\]
\noindent $\textnormal{If~} y_i \in \Omega_{\mathfrak{s}}\cup\Omega_{w}$, by Lemma \ref{lemma1.6}, label the pendant vertices adjacent to $y_i$ in such a way that their label sum is $(1-\deg_s(x))g$.\\
    Suppose  $z\in\Omega_{\mathfrak{p}} $. Then, \[w_{\scriptscriptstyle{\Omega_{\mathfrak{p}}}}(z)  =g.\]
    Suppose  $y_i\in\Omega_{\mathfrak{s}} \cup\Omega_{w}$. Then, 
\begin{align*}
w_{\scriptscriptstyle{\Omega_{\mathfrak{s}} \cup\Omega_{w}}}(y_i) & = \ell_{\scriptscriptstyle {v \in N_{(\Omega_{\mathfrak{s}} \cup \Omega_{w})}(y_i)}}(v)+ \ell_{\scriptscriptstyle{v \in N_{\Omega_{\mathfrak{n}}}(y_i)}}(v)+ \ell_{\scriptscriptstyle{v \in N_{\Omega_{\mathfrak{p}}}(y_i)}}(v) \\
 &= (g+g+\cdots + g)+ 0 + (1-\deg_s(y_i))g\\
 &= \deg_s(y_i)g+ (1-\deg_s(y_i))g\\
 w_{\scriptscriptstyle{\Omega_{\mathfrak{s}} \cup\Omega_{w}}}(y_i) & =g~~\textnormal{for all } y_i \in  \Omega_{\mathfrak{s}}\cup \Omega_{w}.
\end{align*}
Suppose  $x\in\Omega_{\mathfrak{n}}$. Then, 
\begin{align*}
w_{\scriptscriptstyle{\Omega_{\mathfrak{s}} \cup\Omega_{w}}}(x) & = \deg (x)g \\
 &= g  ~~\textnormal{for all } x \in  \Omega_{\mathfrak{n}}.
 \end{align*}
 Therefore, $G$ is $\Gamma$ vertex magic graph with magic constant $\mu =g$.
\end{proof}

Now the following theorems are generalised results of the theorems in \cite{MR4515472}.

\begin{theorem}{\label{thm3}}
Let $\Gamma$ be an Abelian group and $G$ be a graph with all vertices are neutral. If $\ell :V(G)\rightarrow\Gamma$ is a $\Gamma$-vertex magic graph labeling with magic constant $g$  and $G$ contains an edge $uv$ such that $\ell(u)=\ell(v)\neq g$. Then ${G^\dagger}$ is a $\Gamma$-vertex magic graph, where ${G^\dagger}$ is the graph obtained from $G$ by subdividing the edge $uv$ into  $n\equiv 0 \pmod{4}$ times.
\end{theorem}
\begin{proof}
Given that $\ell :V(G) \rightarrow \Gamma$ is a $\Gamma$- vertex magic graph with  magic constant $g$.
 Let $P = (u, u_{1}, u_{2}, u_{3}, \ldots, u_{n-1}, u_{n}, v)$ be a $uv$ path in ${G^\dagger}$. Then $V(G^\dagger)=V(G)\cup \{u_{1}, u_{2}, \ldots, u_{n}\}$.
\noindent Consider a mapping $\ell^{\dagger}:V(G^\dagger)\rightarrow \Gamma$ by, 
%\noindent { Case 1.} $n \equiv 0 \pmod 4$.
\[\ell^\dagger(x)=
\begin{cases} \ell(x) & \textnormal{if}~ x\in V(G)\\
\ell(u)&\textnormal{if}~ x=u_{i},  i \equiv 0 \textnormal{~or~} 1 \pmod{4}\\
%\ell(y) & \textnormal{if}~ v=u_{i},  i \equiv 1\pmod{4} \\
g- \ell(u)&\textnormal{if}~ x=u_i, i \equiv 2 \textnormal{~or~} 3 \pmod{4}.\\

\end{cases}
\]
Then for any $x$ of $G$, $w_{G^\dagger} (x) =g$, and  for any $x$ of the path $P$, $w_{G^\dagger} (x) =\ell(u)+g-\ell(u)=g$.
    \noindent  Clearly, $\ell^\dagger$ is an $\Gamma$-vertex magic labeling of $G^\dagger$ with the magic constant $g$.
\end{proof}
 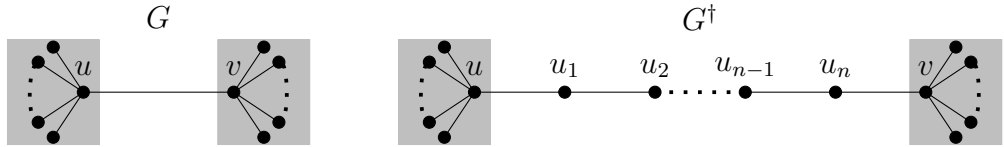
\begin{figure}[htbp]
 	\centering
 	\begin{tikzpicture}[scale=0.4]
 	%\filldraw[fill = black!10, draw = black!10] (1.5,0) circle (4 cm);
 	 
          \filldraw[lightgray,line width=35pt,rounded corners=4pt] (-4,-2)--(-4,2)--cycle;
          \filldraw[lightgray,line width=35pt,rounded corners=4pt] (3,-2)--(3,2)--cycle;
          \filldraw[lightgray,line width=35pt,rounded corners=4pt] (9,-2)--(9,2)--cycle;
          \filldraw[lightgray,line width=35pt,rounded corners=4pt] (26,-2)--(26,2)--cycle;
        
 	\filldraw[fill = black, draw = black] (-3,0) circle (0.2 cm);
 	\filldraw[fill = black, draw = black] (2,0) circle (0.2 cm);
 	\draw (-3,0)--(2,0);
 %	\filldraw[fill = black, draw = black] (-1,0) circle (0.2 cm);
 	\filldraw[fill = black, draw = black] (-4.5,-1) circle (0.2 cm);
 	\filldraw[fill = black, draw = black] (-4.5,1) circle (0.2 cm);
 		\draw [line width=0.5mm, loosely dotted] (-4.5,-1) to[out=180,in=200,looseness=0.5] (-4.5,1);
			\draw [line width=0.5mm, loosely dotted] (3.5,-1) to [out=340,in=30,looseness=0.5] (3.5,1);
 
 	\filldraw[fill = black, draw = black] (3.5,-1) circle (0.2 cm);
 	\filldraw[fill = black, draw = black] (3.5,1) circle (0.2 cm);
 		\filldraw[fill = black, draw = black] (-4,-1.5) circle (0.2 cm);
 			\filldraw[fill = black, draw = black] (-4,1.5) circle (0.2 cm);
 				\draw (-3,0)--(-4,1.5);
 	\draw (-3,0)--(-4,-1.5);
 	\draw (-3,0)--(-4.5,-1);
 	\draw (-3,0)--(-4.5,1);
 	\draw (2,0)--(3.5,-1);
 	\draw (2,0)--(3.5,1);
\filldraw[fill = black, draw = black] (9,1.5) circle (0.2 cm);
 				\draw (10,0)--(9,1.5);
 	\filldraw[fill = black, draw = black] (9,-1.5) circle (0.2 cm);
 				\draw (10,0)--(9,-1.5);
 				\filldraw[fill = black, draw = black] (26,1.5) circle (0.2 cm);
 				\draw (25,0)--(26,1.5);
 	\filldraw[fill = black, draw = black] (26,-1.5) circle (0.2 cm);
 				\draw (25,0)--(26,-1.5);
 					\filldraw[fill = black, draw = black] (3,1.5) circle (0.2 cm);
 				\draw (2,0)--(3,1.5);
 	\filldraw[fill = black, draw = black] (3,-1.5) circle (0.2 cm);
 				\draw (2,0)--(3,-1.5);
 	
\filldraw[fill = black, draw = black] (10,0) circle (0.2 cm);
 \filldraw[fill = black, draw = black] (25,0) circle (0.2 cm);
\filldraw[fill = black, draw = black] (13,0) circle (0.2 cm);
 \filldraw[fill = black, draw = black] (16,0) circle (0.2 cm);
 \filldraw[fill = black, draw = black] (19,0) circle (0.2 cm);
 \filldraw[fill = black, draw = black] (22,0) circle (0.2 cm);
 \filldraw[fill = black, draw = black] (8.5,-1) circle (0.2 cm);
 \filldraw[fill = black, draw = black] (8.5,1) circle (0.2 cm);
 \filldraw[fill = black, draw = black] (26.5,1) circle (0.2 cm);
 \filldraw[fill = black, draw = black] (26.5,-1) circle (0.2 cm);
 \draw[line width=0.5mm,loosely dotted] (16,0)--(19,0);	
 \draw(10,0) -- (16,0);
  \draw(19,0) -- (25,0);

\draw (10,0)--(8.5,-1);
\draw (25,0)--(26.5,1) ;
\draw (25,0)--(26.5,-1) ;
 	\draw (10,0)--(8.5,1);
 		\draw [line width=0.5mm, loosely dotted](8.5,-1) to[out=180,in=200,looseness=0.5] (8.5,1);
			\draw [line width=0.5mm, loosely dotted] (26.5,-1) to [out=340,in=30,looseness=0.5] (26.5,1);
 		
 	\node at (-0.5,2.5) {$G$};
 	\node at (17.5,2.5) { $G^\dagger$};
 	\node at (10,0.8) { $u$};
 	\node at (13,0.8) { $u_1$};
  \node at (16,0.8) { $u_2$};
  \node at (19,0.8) { $u_{n-1}$};
  \node at (22,0.8) { $u_n$};
 	\node at (25,0.8) { $v$};
 		\node at (-3,0.8) { $u$};
 			\node at (2,0.8) { $v$};
 	\end{tikzpicture}
 	\caption{Graphs $G$ and $G^\dagger$}
 \end{figure}

%-------------------
\begin{corollary}
Let $G$ is a $\Gamma$-vertex magic graph with constant $g$. If $G$ contains an edge $uv$ with $\ell(u)=\ell(v)$ and $2\ell(u)=g$, then $G^\dagger$ is also a $\Gamma$-vertex magic, where $G^\dagger$ is the graph obtained from $G$ by subdividing the edge $uv$ for $n$ times.
\end{corollary}
\begin{proof}
Let $G$ be a graph with the edge $uv$. Since $uv$ is subdivided $n$ times, we have a new path $P=(u,u_1,u_2,\ldots,u_n,v)$ connecting $u$ and $v$ in $G^\dagger$. 
%Let  $V(G^\dagger)=V(G)\cup \{u_{1}, u_{2}, \cdots, u_{n}\}$ and  $E(G^\dagger)= E(G)\setminus (xy)  \cup \{(xu_{1}),(u_1u_2), \cdots, (u_{n}y)\}$.
\noindent Consider a mapping $\ell^{\dagger}:V(G^\dagger)\rightarrow \Gamma$ given by, 
\[\ell^\dagger(x)=
\begin{cases} \ell(x) & \textnormal{if}~ x\in V(G)\\
\ell(u)&\textnormal{if}~ x\in V(P)\setminus\{u,v\}.
\end{cases}
\]
\noindent  It is not hard to verify that $\ell^\dagger$ is an $\Gamma$-vertex magic labeling of $G^\dagger$ with the magic constant $g$.
\end{proof}
%------------------------

\begin{theorem}
 Let $\Gamma$ be any Abelian group with $|\Gamma|\geq 3$ and let $g \in \Gamma\setminus\{0\}$. Let $t_i\ge 2$. Suppose $G$ is $\Gamma$-vertex magic graph with labeling $\ell$ and magic constant $g$. If there exists an edge $uv$ in $G$ with $\ell (u)=\ell(v)=g$, then the graph $G^\dagger$ obtained from $G$ by subdividing the edge $uv$ into $n$ times and by attaching $t_i$ pendant vertices at each new vertices $x_i$, is an $\Gamma$-vertex magic graph with the same magic constant $g$. \end{theorem}
 \begin{proof}
 Define $\ell^\dagger: V(G^\dagger)\rightarrow \Gamma$ by $\ell^\dagger(x_i)=g$ and label the $t$ pendant vertices adjacent to each $x_i$ such that its label sum is $-g$ and $\ell^\dagger(w)=\ell(w)$ for all $w
 \in V(G)$. Clearly $\ell^\dagger$ is an $\Gamma$-vertex magic of $G^\dagger$ with magic constant $g$.
 \end{proof}
 \begin{figure}[htbp]
 	\centering
 	\begin{tikzpicture}[scale=0.4]
 	%\filldraw[fill = black!10, draw = black!10] (1.5,0) circle (4 cm);
 	 
          \filldraw[lightgray,line width=35pt,rounded corners=4pt] (-4,-2)--(-4,2)--cycle;
          \filldraw[lightgray,line width=35pt,rounded corners=4pt] (3,-2)--(3,2)--cycle;
          \filldraw[lightgray,line width=35pt,rounded corners=4pt] (9,-2)--(9,2)--cycle;
          \filldraw[lightgray,line width=35pt,rounded corners=4pt] (26,-2)--(26,2)--cycle;
       
 	\filldraw[fill = black, draw = black] (-3,0) circle (0.2 cm);
 	\filldraw[fill = black, draw = black] (2,0) circle (0.2 cm);
 	\draw (-3,0)--(2,0);
 %	\filldraw[fill = black, draw = black] (-1,0) circle (0.2 cm);
 	\filldraw[fill = black, draw = black] (-4.5,-1) circle (0.2 cm);
 	\filldraw[fill = black, draw = black] (-4.5,1) circle (0.2 cm);
 		\draw [line width=0.5mm, loosely dotted] (-4.5,-1) to[out=180,in=200,looseness=0.5] (-4.5,1);
			\draw [line width=0.5mm, loosely dotted] (3.5,-1) to [out=340,in=30,looseness=0.5] (3.5,1);
 
 	\filldraw[fill = black, draw = black] (3.5,-1) circle (0.2 cm);
 	\filldraw[fill = black, draw = black] (3.5,1) circle (0.2 cm);
 		\filldraw[fill = black, draw = black] (-4,-1.5) circle (0.2 cm);
 			\filldraw[fill = black, draw = black] (-4,1.5) circle (0.2 cm);
 				\draw (-3,0)--(-4,1.5);
 	\draw (-3,0)--(-4,-1.5);
 	\draw (-3,0)--(-4.5,-1);
 	\draw (-3,0)--(-4.5,1);
 	\draw (2,0)--(3.5,-1);
 	\draw (2,0)--(3.5,1);
\filldraw[fill = black, draw = black] (9,1.5) circle (0.2 cm);
 				\draw (10,0)--(9,1.5);
 	\filldraw[fill = black, draw = black] (9,-1.5) circle (0.2 cm);
 				\draw (10,0)--(9,-1.5);
 				\filldraw[fill = black, draw = black] (26,1.5) circle (0.2 cm);
 				\draw (25,0)--(26,1.5);
 	\filldraw[fill = black, draw = black] (26,-1.5) circle (0.2 cm);
 				\draw (25,0)--(26,-1.5);
 					\filldraw[fill = black, draw = black] (3,1.5) circle (0.2 cm);
 				\draw (2,0)--(3,1.5);
 	\filldraw[fill = black, draw = black] (3,-1.5) circle (0.2 cm);
 				\draw (2,0)--(3,-1.5);
 	
\filldraw[fill = black, draw = black] (10,0) circle (0.2 cm);
 \filldraw[fill = black, draw = black] (25,0) circle (0.2 cm);
\filldraw[fill = black, draw = black] (13,0) circle (0.2 cm);
 \filldraw[fill = black, draw = black] (16,0) circle (0.2 cm);
 \filldraw[fill = black, draw = black] (19,0) circle (0.2 cm);
 \filldraw[fill = black, draw = black] (22,0) circle (0.2 cm);
 \filldraw[fill = black, draw = black] (8.5,-1) circle (0.2 cm);
 \filldraw[fill = black, draw = black] (8.5,1) circle (0.2 cm);
 \filldraw[fill = black, draw = black] (26.5,1) circle (0.2 cm);
 \filldraw[fill = black, draw = black] (26.5,-1) circle (0.2 cm);
% \draw[line width=mm,loosely dotted] (16,0)--(19,0);
 \draw(10,0) -- (16,0);
  \draw(19,0) -- (25,0);

\draw (10,0)--(8.5,-1);
\draw (25,0)--(26.5,1) ;
\draw (25,0)--(26.5,-1) ;
 	\draw (10,0)--(8.5,1);
 		\draw [line width=0.5mm, loosely dotted](8.5,-1) to[out=180,in=200,looseness=0.5] (8.5,1);
			\draw [line width=0.5mm, loosely dotted] (26.5,-1) to [out=340,in=30,looseness=0.5] (26.5,1);
 		
 	\node at (-0.5,2.5) {$G$};
 	\node at (17.5,3) { $G^\dagger$};
 	\node at (10,0.8) { $u$};
 	\node at (13,0.8) { $x_1$};
  \node at (16,-0.8) { $x_2$};
  \node at (19,0.8) { $x_{n-1}$};
  \node at (22,-0.8) { $x_n$};
  \filldraw[fill = black, draw = black] (11.5,-2) circle (0.2 cm);
  \draw (13,0)--(11.5,-2);
    \draw (13,0)--(12.5,-2);
      \draw (13,0)--(14.5,-2);
      \node at (13.5,-2){$\cdots$};
       \node at (16.5,2){$\cdots$};
       \node at (19.5,-2){$\cdots$};
        \node at (22.5,2){$\cdots$};
 %     \draw[loosely dotted](12.5,-2)--(14.5,-2);
  \filldraw[fill = black, draw = black] (12.5,-2) circle (0.2 cm);
  \filldraw[fill = black, draw = black] (14.5,-2) circle (0.2 cm);
    \draw (16,0)--(14.5,2);
    \draw (16,0)--(15.5,2);
    \draw (16,0)--(17.5,2);
    \draw (19,0)--(17.5,-2);
    \draw (19,0)--(18.5,-2);
    \draw (19,0)--(20.5,-2);
    \draw (22,0)--(20.5,2);
     \draw (22,0)--(21.5,2);
      \draw (22,0)--(23.5,2);
   \filldraw[fill = black, draw = black] (14.5,2) circle (0.2 cm);
  \filldraw[fill = black, draw = black] (15.5,2) circle (0.2 cm);
  \filldraw[fill = black, draw = black] (17.5,2) circle (0.2 cm);
   \filldraw[fill = black, draw = black] (17.5,-2) circle (0.2 cm);
  \filldraw[fill = black, draw = black] (18.5,-2) circle (0.2 cm);
  \filldraw[fill = black, draw = black] (20.5,-2) circle (0.2 cm);

   \filldraw[fill = black, draw = black] (20.5,2) circle (0.2 cm);
  \filldraw[fill = black, draw = black] (21.5,2) circle (0.2 cm);
  \filldraw[fill = black, draw = black] (23.5,2) circle (0.2 cm);
 	\node at (25,0.8) { $v$};
 		\node at (-3,0.8) { $u$};
 			\node at (2,0.8) { $v$};
    \node at (17.5,0) {$\cdots$};
 	\end{tikzpicture}
 	\caption{Graphs $G$ and $G^\dagger$}
 \end{figure}
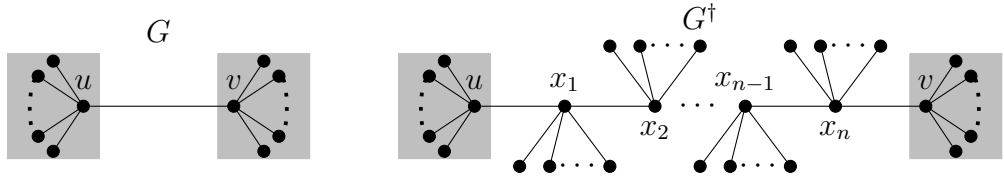
     %----------------
\section{Group vertex magicness of product graphs}
In this section, the group vertex magicness of certain product graphs are determined. We adopt the following notation and terminology. 

\begin{proposition}{\label{note4.1}}
Let $G$ be a graph in which either all vertices are of odd degrees or all vertices are of even degrees and $\Gamma$ an Abelian group with at least one involution. Then $G$ is a $\Gamma$-vertex magic graphs. In addition, if $G$ is regular then $G$ is a group vertex magic.
\end{proposition}
\begin{proof}

Consider the mapping $\ell$ from $V(G)$ to $\Gamma$ such that $\ell(x)=g$, for all $x \in V(G)$.  Suppose $g$ is an involution element in $\Gamma$. \\
Case 1. If $G$ is a graph with all vertices are even degrees, then, \[w_G(x)=\underbrace{g+g+\cdots+g}_{ \textnormal{(even number of times)~}} = 0, \textnormal{~for all~} x \in V(G).\] 
Case 2. If $G$ is a graph with all vertices are odd degrees, then, \[w_G(x)=\underbrace{g+g+\cdots+g}_{\textnormal{(odd number of times)~}} = g, \textnormal{~for all~} x \in V(G).\] 
In both cases, $G$ is  $\Gamma$-vertex magic.
\end{proof}
%--------
 \begin{observation}
      Let $\Gamma$ be a Abelian group. If $\ell$ be a $\Gamma$-vertex magic labeling of $G$ with magic constant $\mu$,  then there exits a $\Gamma$-vertex magic labeling $\ell^\dagger$ of $G$ with magic constant $-\mu$.
  \end{observation}
  \begin{proof}
      Given that $\ell$ be a $\Gamma$-vertex magic labeling of $G$. Then for any $u\in V(G)$, \[\sum_{v\in N(u)}\ell(u)=\mu\]
    Define a labeling $\ell^{\dagger}:V(G)\rightarrow \Gamma$ such  that $\ell^\dagger(u)=-\ell(u) \textnormal{~for all~} u \in V(G)$. Then \[\sum_{v\in N(u)}\ell^\dagger(u)=\sum_{v\in N(u)}-\ell(u)=-\sum_{v\in N(u)}\ell(u)=-\mu\]
  \end{proof}
  
\begin{observation} Let $G$  be a graph.
   If $G$ is a group vertex magic graph with all vertex labels are same if and only if $G$ is regular.
\end{observation}
\begin{proof}
    Let $g$ be a non identity element in an Abelian group.  Assume that $G$ is a group vertex magic graph with all vertex labels are $g$. Then, for each vertex $v_i\in G$, $w_G(v_i)=\deg(v_i)g$. Therefore, the degree must be same for all the vertices, since $G$ is a group vertex magic.\\
    The converse, is from Proposition \ref{note4.1}.
\end{proof}

% \begin{proposition}Let $G$ be $r$-regular graph and $M$ a perfect matching in $G$. Then $G\setminus M$ is a group vertex magic graph.  \textcolor{red}{note G-M is another regular graph.}
% \end{proposition}
% \begin{proposition} \textcolor{red} {(true? check! think of all labels of v are labelled g)Let $G$ be any graph and $M$ a perfect matching in $G$. If G is group vertex magic with magic constant $kg$, then $G\setminus M$ is also a group vertex magic graph with magic constant $(k-1)$g.  } 
% \end{proposition}
% \begin{proof}
% From Proposition {\ref{note4.1}}, $G$ is a group vertex magic graph with $\mu=rg$. Let $M$ be any perfect matching in $G$. Then degree of each vertex of $G \setminus M$ is reduced by 1. Suppose, $G\setminus M$ preserve the same labeling as in $G$. Then $G$ is a group vertex magic graph with $\mu=(r-1)g$.
% \end{proof}
%-------------

%-------------
The following theorems deal with the graph products of a graph $G$ with some well known classes of graphs.
%---------------
\begin{theorem}
    Let $G$ and $G^\dagger$ be two graphs. If $G$ is a $\Gamma$-vertex magic graph with labeling $\ell$ and magic constant $\mu=g$ such that $\ell(x)=g$ for all $x\in V(G)$. Also, $G^\dagger$ is a  $\Gamma^\dagger$-vertex magic graph with labeling $\ell^\dagger$ and magic constant $\mu=g^\dagger$ such that $\Sigma_{\ell^\dagger(G^\dagger)}=g^\dagger$. Then $G\odot G^\dagger$ is a $\Gamma\times \Gamma^\dagger$- vertex magic graph.
\end{theorem}
   \begin{proof}
     Let  $G^\ddagger= G\odot G^\dagger$. Consider the mapping $\ell^\ddagger: V(G^\ddagger)\rightarrow \Gamma\times \Gamma^\dagger$ such that,\[\ell^{\ddagger}(x)=
     \begin{cases}
     (\ell(x),0) & \textnormal{If~~} x \in V(G)\\
     (0, \ell^\dagger(x)) & \textnormal{If~~} x \in V(G^\dagger).
     \end{cases}\]
     Suppose $x \in V(G)$. Then, $w_{G^\ddagger}(x)=(\sum_{y\in N_{G}(x)} \ell(y),0)+(0,\Sigma_{\ell^\dagger(G^\dagger)})=(g,0)+(0,g^\dagger)=(g,g^\dagger)$.\\
     Suppose $x \in V(G^\dagger)$. Then, $w_{G^\ddagger}(x)=(\ell(x),0)+(0,\sum_{y\in N_{G^\dagger}(x)} \ell(y))=(g,0)+(0,g^\dagger)=(g, g ^\dagger).$ 
     Hence, $G^\ddagger$ is a $\Gamma\times \Gamma^\dagger$-vertex magic graph.
   \end{proof} 
%---------------
\begin{theorem}
Let $\Gamma$ be any Abelian group with $|\Gamma|\geq 3$ and $G$ be a $r$-regular graph. Then $G\odot K_m$ is a group vertex magic graphs. 
\end{theorem}
\begin{proof}
Let $|V(G)|=n$ and  $G^\dagger= G\odot K_m$. Suppose $G^\dagger$ is a $\Gamma$ vertex magic labeling then its must have a labeling $\ell$ such that $w_{G^\dagger}(x)=g , g \in \Gamma \textnormal{ for all } x\in V(G^\dagger)$. Assume that $\ell(x) = g_{1}, g_{1}\in \Gamma $ whenever $x\in G \subseteq G^\dagger$ and $\ell(x) = g_{2}, g_2 \in \Gamma $ whenever $x\in K_m \subseteq G^\dagger$. Then $w_{G^\dagger}(x)= rg_1 + mg_2$ and when $x\in G$ and $w_{G^\dagger}(x)= g_1 + (m-1)g_2$ and when $x\in K_m$. Therefore,
\[rg_1+mg_2 = g_1 +(m-1)g_2 \textnormal{~and hence~} g_2 =-(r-1)g_1.\] Then \[\ell(x) =
\begin{cases}
g_1 & \textnormal{If~} x \in V(G)\\
 -(r-1)g_1 & \textnormal{If~} x \in V(K_m).
\end{cases}\]
Thus $w_{G^\dagger}(x)=g_1(r-m(r-1))$. Hence $G\odot K_m$ is a $\Gamma$- vertex magic graph. Also the given labeling satisfied for all additive Abelian group with $|\Gamma|\geq 3$. Therefore $G\odot K_m$ is a group vertex magic graph. 
\end{proof}
%------------
\begin{observation}
    Let $\Gamma$ be Abelian group with $o(\Gamma)\ge 3$ and $G$ be a $r$-regular graph with order $n$. If $g\in \Gamma$ such that $n \equiv 1 \pmod{o(g)}$, then $G \odot G$ is $\Gamma$-vertex magic graph.
\end{observation}
\begin{proof}
If we label all the vertices of $G\odot G$ by a non-zero element $g$ of $\Gamma$, then the weight of any vertex is $rg + ng = rg+g =(r+1)g$.
\end{proof}
%---------
% \begin{observation}
%    \color{red}{(true? think) Let $\Gamma$ be Abelian group with $o(\Gamma)\ge 3$ and $G$ be a $r$-regular graph with order $n$. \st{If $g\in \Gamma$ such that $n \equiv 1 \pmod{o(g)}$, then} $G \odot G$ is $\Gamma$-vertex magic graph.}
% \end{observation}
% \begin{proof}
% If we label all the vertices of $G\odot G$ by a non-zero element $g$ of $\Gamma$, then the weight of any vertex is $(r+n)g$.
% \end{proof}
%----------------
% \begin{observation}
%     Let $G$ be a $r$-regular graph of order $m$ and $H_1, H_2, \cdots, H_m$ have regularities $k$ and order $n$. Let $g$ be a non identity element in $\Gamma$ such that $(n+r)\equiv (k+1) \pmod{o(g)}$. Then $G\tilde{\circ} \bigl(\mathop{\Lambda}\limits^m_{i=1}H_i\bigr)$ is $\Gamma$-vertex magic graph.
% \end{observation}
%---------------
\begin{observation}
Let $G$ be  $r$-regular graph of order $m$, and $H_1, H_2, \cdots, H_m$ are $r_1, r_2, \cdots, r_m$ regularity respectively, where each $H_i$ is of order $n$. If there is an non-identity element $g$  in $\Gamma$ such that $(r_i+n) \equiv (r_j+1) \pmod{(o(g)}$. Then $G\tilde{\circ} \bigl(\mathop{\Lambda}\limits^m_{i=1}H_i\bigr)$ is $\Gamma$-vertex magic graph.
\end{observation}

\begin{theorem}
Let $G$ be a graph with $m$ vertices and $H$ be an $r$-regular graph with $n$ vertices. Let $p$ be any prime  such that  $n \equiv 0 \pmod p$. If $\Gamma$ is an Abelian group with  $p^2$ elements, then $G[H]$ is $\Gamma$-vertex magic.
\end{theorem}
\begin{proof}
\noindent { {Case 1.}} Suppose $G$ is a  $k$ -regular graph.\\
Given that $H$ is a $r$- regular graph with $n$ vertices, then $G[H]$ is $r +k m$ regular graph. Now,  by Proposition \ref{note4.1},  there exists an Abelian group $\Gamma$ such that $\ell :V(G[H])\rightarrow\Gamma$ is a group vertex magic graph labeling with  $\mu=g$, where $g \in \Gamma$. \\
 \noindent { {Case 2.}} Suppose $G$ is a  non-regular graph.\\
 \noindent Since $H$ is an $r$-regular graph, by Proposition \ref{note4.1},  $H$ is group vertex magic with a magic constant $rg$, for some non-zero element $g, g \in \Gamma$.\\
  \noindent Suppose  $\ell$ is a mapping from the vertex set of $G[H]$ to $\Gamma$. \\ Let  $H_{1}, H_{2}, \cdots, H_{m}$ be the $m$ copies of $H$ in  $G[H]$.  Let $\ell(x)=g$ for all $x \in H$.  
  Then for each $H_{i}$,  $ \Sigma_{\ell(H_i)}=0,$ because $n \equiv 0 \pmod p$. \\
Thus $G[H]$ is a group vertex magic graph with  magic constant $rg$ 
\end{proof}
%---------------------------
\begin{theorem} {\label{thm3.7}}
   Let $\Gamma$ be any Abelian group with $|\Gamma|\ge 3$. Let $G$ be an $r$-regular graph with $n$ vertices. If $H_{1}, H_{2}, \cdots, H_{n}$ are $\Gamma$-vertex magic  graphs with $\mu(H_i)=g_1$ for all $i$,  and  $\Sigma_{\ell_i(H_i)}=g$ for all $i$, where $g, g_1 \in \Gamma$. Then $G[H_{1}, H_{2}, \cdots, H_{n}]$ is $\Gamma$-vertex magic graph.
\end{theorem}
%%%%%%%%%%%%%%%%%%%%%% TYPOS -------------
\begin{proof} 
Let $G^\dagger=G[H_1, H_2, \cdots, H_n]$. Since $G$ is regular, by Proposition \ref{note4.1}, $G$ is a group vertex magic graph. Suppose $V(G)=\{v_1, v_2, \cdots, v_n\}$. Then $V(G^{\dagger}$ is defined as $\{v_i^j, 1 \le i \le n , 1 \le j \le |V(H_i)|\}$.\\ Let $\ell_i$ be the corresponding $\Gamma$-vertex magic labeling of $H_i, 1 \le i \le n$ with $\mu(H_i)=g_1$ for all $i$  and  $\Sigma_{\ell_i(H_i)}=g$ for all $i$, where $g\in \Gamma$.\\
Now, consider a mapping  $\ell^\dagger$ from the vertex set of $G^\dagger$ to $\Gamma$ such that each $H_i$ preserve the same labeling as in $\ell_i$ for any $i$. That is, 
$\ell^\dagger(x_i^j)=\ell_{H_j}(x_i)$ for all $x_i^j \in V(G^\dagger)$  
Then for any vertex  $x \in G^\dagger$ , $w_{G^\dagger}(x)= g_1 + rg$. Therefore $G^\dagger$ is a  $\Gamma$ vertex magic graph.
\end{proof}
%--------------
\begin{observation}
Let $G$ be a any graph with $n$ vertices.  If $H_{1}, H_{2}, \cdots, H_{n}$ are $\Gamma$-vertex magic  graphs with $\mu(H_i)=g$ for all $i$,  and  $\Sigma_{\ell_i(H_i)}=0$ for all $i$. Then $G[H_{1}, H_{2}, \ldots, H_{n}]$ is $\Gamma$-vertex magic graph.
\end{observation}
%-----------------------
\begin{theorem}
Let $G$ be a graph  with $n$ vertices and   $H_{1}, H_{2},\ldots, H_{n}$ be any graphs.  Let $A$ be an Abelian group and $g\in A$. Suppose $H_{1}, H_{2}, \cdots, H_{n}$ are $\Gamma$-vertex magic  graphs with $\mu(H_i)=g$ for all $i$  and  $\Sigma_{\ell_i(H_i)}=g$ for all $i$. Then $G[H_{1}, H_{2}, \ldots, H_{n}]$ is a $\Gamma$-vertex magic graph with magic constant $g$ if and only if degree of all the vertices of  $G[H_{1}, H_{2}, \ldots, H_{n}]$ is an integral multiple of $o(g)$.
\end{theorem}
\begin{proof} Take $G^\dagger=G[H_{1}, H_{2}, \ldots, H_{n}]$.
Given that $H_{1}, H_{2}, \ldots, H_{n}$ are $\Gamma$-vertex magic graphs with magic constant $g$. Suppose $\ell_1, \ell_2, \ldots, \ell_n$ are the corresponding $\Gamma$-vertex magic labelings of $H_{1}, H_{2}, \ldots, H_{n}$. Then for any vertex $x\in H_{i}$,
\begin{align}{\label{eqn3}}
  w_{\scriptscriptstyle{H_i}}(x)=  g=\sum_{\scriptscriptstyle{y\in N_{H_i}(x)}}\ell_{i}(y)
\end{align}
Assume that $\ell^\dagger:V(G^\dagger) \rightarrow \Gamma$,  is a $\Gamma$-vertex magic graph with magic constant $g$. Also, $\ell^\dagger(x_i^j)=\ell_{H_j}(x_i)$ for all $x_i^j \in V(G^\dagger)$.\\
Case 1. Suppose the graph $G$ is regular. Then by Theorem \ref{thm3.7}, $G^\dagger$ is $\Gamma$-vertex magic graph. Hence  $w_{\scriptscriptstyle{G^\dagger}}(x)= (r+1)g$, for all $x\in V(G^\dagger)$. Thus $g=rg+g$, implies $rg=0$. Therefore $r\equiv 0 \pmod{o(g)}$

\noindent Case 2. Suppose the  graph $G$ is not regular. Let $x\in V( G^\dagger)$. If $x\in V(H_i)$,  then,
\begin{align*}
   w_{\scriptscriptstyle{G^\dagger}}(x)= g & = \sum_{y\in N_{H_i}(x)}\ell^\dagger(x)+\sum_{H_j\in N(H_i)}(\Sigma_{\ell^\dagger(H_j)})\\
   g & = \sum_{\scriptscriptstyle{x\in N_{H_1}(x_{1})}}\ell_{1}(x) + \sum_{H_j\in N(H_i)}(\Sigma_{\ell_{j}(H_j)})\\
  g & = g +\deg(x)g \\
 0 & = \deg(x)g.
\end{align*}
Therefore, it is clear that, $\deg(x)\equiv 0 \pmod{o(g)}$.
\\ For the converse, consider $\ell^\dagger:V(G^\dagger)\rightarrow \Gamma$ such that the  mapping $\ell^\dagger$ sustain the same labeling $\ell_{i}$ in each $H_i,  1\le i \le n$. \\
Let $x\in V( G^\dagger)$, in particular $x\in V(H_{i})$. Then,
\begin{align*}{\label{eqn4}}
   w_{\scriptscriptstyle{G^\dagger}}(x) & = \sum_{x\in N_{H_i}(x)}\ell^\dagger(x)+\sum_{H_j\in N(H_i)}(\Sigma_{\ell^\dagger(H_j)})\\
   & = \sum_{\scriptscriptstyle{x\in N_{H_i}(x)}}\ell_{i}(x) + \sum_{H_j\in N(H_i)}(\Sigma_{\ell_{j}(H_j)})\\
  w_{\scriptscriptstyle{G^\dagger}}(x)& = g +\deg(x)g,  \textnormal{~Since~} \deg\equiv 0\pmod{o(g)}   \\
 w_{\scriptscriptstyle{G^\dagger}}(x)& = g. 
\end{align*}
Hence for any vertices in $G^\dagger$, we get the same weight. Therefore $G^\dagger$ is a $\Gamma$-vertex magic graph.
\end{proof}
%------------------------
\begin{theorem}
Let $\Gamma$ and $\Gamma^\dagger$ be two Abelian groups. Let $G$ and $G^{\dagger}$ be two graphs of order $m$ and 
$n$, respectively. If $G$ is a $\Gamma$ vertex magic graph with magic constant $g$ and 
$G^{\dagger}$ is a $\Gamma^\dagger$-vertex magic graph with magic constant $g^{\dagger}$, then $G \square G^{\dagger}$ is a  $\Gamma \times\Gamma^\dagger$-vertex magic graph with magic constant $(tg, tg^{\dagger})$, where $t=m+n$
.\end{theorem}
\begin{proof}
Given that $\ell:V(G)\rightarrow \Gamma$ is a $\Gamma$-vertex magic graph with magic constant $g$. Then,
\begin{equation}{\label{eqn1}}
    \sum_{x \in V(G)} \deg(x)\ell(x)= mg. 
\end{equation}
Again, $\ell^{\dagger}:V(G^{\dagger})\rightarrow \Gamma^\dagger$ is a $\Gamma^\dagger$-vertex magic graph with magic constant $g^{\dagger}$. Then,
\begin{equation}{\label{eqn2}}
    \sum_{y \in V(G^{\dagger})} \deg(y)\ell^{\dagger}(y)= ng^{\dagger}.
\end{equation}
Define a mapping $\ell^{\ddagger}$ from the vertex set of $G \square G^{\dagger}$ to $\Gamma \times\Gamma^\dagger$ by
\[\ell^{\ddagger}(x,y)=(\ell(x), \ell^{\dagger}(y)), x\in V(G), y \in V(G^{\dagger}).\]
Then the weight of any vertex $(x,y)$ \[ w^{\ddagger}(x,y) = (ng + \sum_{x \in V(G)} \deg(x)\ell(x), mg^\dagger +\sum_{y \in V(G^{\dagger})} \deg(y)\ell^{\dagger}(y))\]
From equations $\ref{eqn1}$ and $\ref{eqn2}$, $w^{\ddagger}(x,y)=((m+n)g, (m+n)g^{\dagger})$.
\end{proof}

\begin{theorem}
If $G$ is the graph defined in Theorem ${\ref{thm2.2}}$ and $H$ is any graph with at least one pendant vertex, then the direct product  $G\times H$ is a not group vertex magic graph.

\end{theorem}
\begin{proof}
Let $G$ and $H$ be two graph with the vertices $u_{1}, u_{2}, \cdots, u_{m}$  and  $v_{1}, v_{2}, \cdots, v_{n}$, respectively. Suppose $u_{1}$ is the pendant vertex in $G$ and $u_{3}$ is of degree 2 in  $G$. Let $v_{1}$ is the pendant vertex in $H$. Then $(u_1,v_1)$ is the pendant vertex in  $G\times H$ and  $(u_3,v_1) \in N(N(u_{1}v_{1}))$ is the vertex with degree $2$ in $G\times H$. By Theorem ${\ref{thm2.2}}$, $G\times H$ is a not group vertex magic.
\end{proof}

\begin{example}
$P_{5}\times P_{3}$ is not group vertex magic.
\end{example}
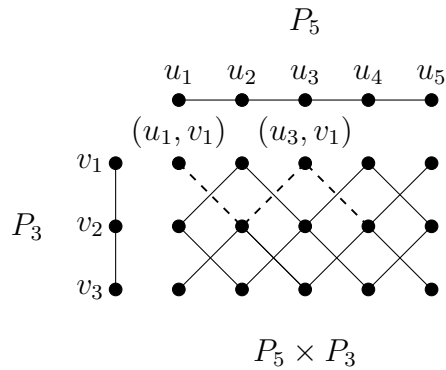
\begin{figure}[htbp]
	\centering
	\begin{tikzpicture}[scale=0.42]
		\filldraw[fill=black](-6,4) circle (0.2cm);
		\filldraw[fill=black](-6,2) circle (0.2cm);
			\filldraw[fill=black](-6,0) circle (0.2cm);
   \draw  (-6,4) -- (-6,0);
\filldraw[fill=black](-4,6) circle (0.2cm);
		\filldraw[fill=black](-2,6)  circle (0.2cm);
			\filldraw[fill=black](0,6) circle (0.2cm);
			\filldraw[fill=black](2,6)  circle (0.2cm);
			\filldraw[fill=black](4,6) circle (0.2cm);
   \draw  (-4,6) -- (4,6);
		\filldraw[fill=black](-4,0) circle (0.2cm);
		\filldraw[fill=black](-2,0) circle (0.2cm);
			\filldraw[fill=black](0,0) circle (0.2cm);
				\filldraw[fill=black](2,0) circle (0.2cm);
					\filldraw[fill=black](4,0) circle (0.2cm);
		\filldraw[fill=black](-4,2) circle (0.2cm);
		\filldraw[fill=black](-2,2) circle (0.2cm);
			\filldraw[fill=black](0,2) circle (0.2cm);
				\filldraw[fill=black](2,2) circle (0.2cm);
					\filldraw[fill=black](4,2) circle (0.2cm);
						\filldraw[fill=black](-4,4) circle (0.2cm);
		\filldraw[fill=black](-2,4) circle (0.2cm);
			\filldraw[fill=black](0,4) circle (0.2cm);
				\filldraw[fill=black](2,4) circle (0.2cm);
					\filldraw[fill=black](4,4) circle (0.2cm);		
\draw  (-2,2) -- (0,0);
\draw  (-2,4) -- (2,0);
\draw  (-2,2) -- (0,0);
\draw  (2,2) -- (4,0);
\draw[line width=0.25mm,dashed](2,2) -- (0,4);
\draw  (4,4) -- (0,0);
\draw  (2,4) -- (-2,0);
\draw  (-2,2) -- (-4,0);
\draw  (-2,4) -- (-4,2);
\draw(-2,0) -- (-4,2);
\draw[line width=0.25mm,dashed]  (-4,4) -- (-2,2);
\draw[line width=0.25mm,dashed]  (0,4) -- (-2,2);
\draw  (2,4) -- (4,2);
\draw   (2,0)--  (4,2);

\node at (0, 8.5){$P_{5}$};
\node at (-4, 6.8){$u_1$};
\node at (-4, 5){$(u_1,v_1)$};
\node at (0, 5){$(u_3,v_1)$};
\node at (-2, 6.8){$u_2$};
\node at (0, 6.8){$u_3$};
\node at (2, 6.8){$u_4$};
\node at (4, 6.8){$u_5$};
\node at (-8.8,2 ){$ P_{3}$};
\node at (-6.8,2 ){$ v_{2}$};
\node at (-6.8,4 ){$ v_{1}$};
\node at (-6.8,0 ){$ v_{3}$};
\node at (0,-2 ){$ P_{5} \times P_{3}$};
 %\begin{scope}[on background layer]   
  %        \filldraw[lightgray,line width=10pt,rounded corners=2pt] (-4,3.9)--(-4,4.1)--cycle;
  %       \filldraw[lightgray,line width=15pt,rounded corners=2pt] (0,3.9)--(0,4.1)--cycle;
  %  \end{scope}

	\end{tikzpicture}
	\end{figure}

%----------------------

\section{Group vertex magicness of trees}{\label{sec.3}}

In this section, we prove a characterization of the group vertex magicness of trees of diameter $4$ and $5$, for any infinite Abelian group $\Gamma$. We consider three cases based on the existence of the torsion elements in $\Gamma$.
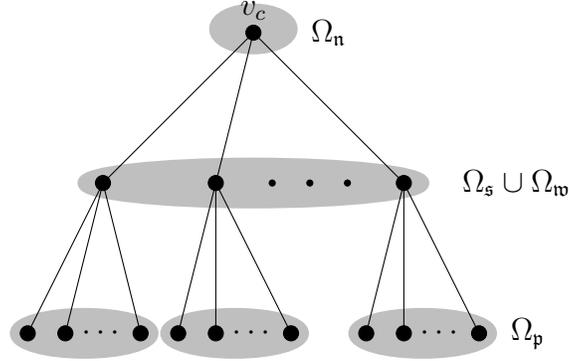
\begin{figure}[htbp]
\centering
   \begin{tikzpicture}[scale=0.5]
   \filldraw[lightgray,line width=5pt,rounded corners=4pt]  (0,0.1) ellipse (1cm and 0.5cm);
    \filldraw[lightgray,line width=5pt,rounded corners=4pt]  (0,-4) ellipse (4.5cm and 0.5cm);
     \filldraw[lightgray,line width=5pt,rounded corners=4pt]  (-4.5,-8) ellipse (1.8cm and 0.5cm);
     \filldraw[lightgray,line width=5pt,rounded corners=4pt]  (-0.5,-8) ellipse (1.8cm and 0.5cm);
     \filldraw[lightgray,line width=5pt,rounded corners=4pt]  (4.5,-8) ellipse (1.8cm and 0.5cm);
    \filldraw[fill=black](0,0) circle (0.2cm);
    \filldraw[fill=black](-4,-4) circle (0.2cm);
    \filldraw[fill=black](-1,-4) circle (0.2cm);
    \filldraw[fill=black](4,-4) circle (0.2cm);
    \filldraw[fill=black](-6,-8) circle (0.2cm);
     \filldraw[fill=black](-5,-8) circle (0.2cm);
      \filldraw[fill=black](-3,-8) circle (0.2cm);
       \filldraw[fill=black](-2,-8) circle (0.2cm);
        \filldraw[fill=black](-1,-8) circle (0.2cm);
         \filldraw[fill=black](1,-8) circle (0.2cm);
          \filldraw[fill=black](6,-8) circle (0.2cm);
           \filldraw[fill=black](3,-8) circle (0.2cm);
            \filldraw[fill=black](4,-8) circle (0.2cm);
            \draw(0,0)--(4,-4);
\draw(0,0)--(-4,-4);
\draw(0,0)--(-1,-4);
\draw(-6,-8)--(-4,-4);
 \draw(-5,-8)--(-4,-4); 
 \draw(-3,-8)--(-4,-4);
 \draw(-2,-8)--(-1,-4);
 \draw(-1,-8)--(-1,-4);
 \draw(1,-8)--(-1,-4);
 \draw(6,-8)--(4,-4);
  \draw(3,-8)--(4,-4);
   \draw(4,-8)--(4,-4);
 \filldraw[fill=black](0.5,-4) circle (0.08cm); 
  \filldraw[fill=black](1.5,-4) circle (0.08cm); 
   \filldraw[fill=black](2.5,-4) circle (0.08cm); 
\node at (-4,-8){$\cdots$};
 \node at (0,-8){$\cdots$}; 
 \node at (5,-8){$\cdots$};
 \node at (0,0.6){$v_c$};
 \node at (2,0){$\Omega_{\mathfrak{n}}$};
 \node at (7,-4){$\Omega_{\mathfrak{s}}\cup \Omega_{\mathfrak{w}}$};
 % \node at (5,-9.2){$\Omega_{\mathfrak{p}}$};
  \node at (7.3,-8){$\Omega_{\mathfrak{p}}$};
  %\node at (-4,-9.2){$\Omega_{\mathfrak{p}}$};
\end{tikzpicture}
\caption{A typical tree with diameter $4$}
\end{figure}

\begin{theorem}Let $T$ be a tree of diameter $4$ with the centre vertex $v_c$ and $\Gamma$ be an infinite Abelian torsion group. Then $T$ is  $\Gamma$-vertex magic if and only if $T$ satisfies one of the criteria below:\\
		\textnormal{(i)} Every non-pendant vertex of $T$ is in $\Omega_{\mathfrak{s}}(T)$ \\ %a strong support vertex.\\
		\textnormal{(ii)} $v_c\in \Omega_{\mathfrak{w}}(T)$, $d(v_c)\not\equiv 2 \pmod {e(\Gamma)}$, and  all other non-pendant vertices are in $\Omega_{\mathfrak{s}}(T)$ \\ %strong support vertices.\\
		\textnormal{(iii)} $v_c \in \Omega_{\mathfrak{n}}(T)$  and $\gcd(d(v_c)-1,e(\Gamma))\neq 1.$
\end{theorem}
\begin{proof}
    Let $T$ be a $\Gamma$-vertex magic  tree with a labeling $\ell$ and a magic constant $g$. There are two cases.\\
    { Case 1}: Suppose $v_c \in (\Omega_{\mathfrak{s}}\cup\Omega_{\mathfrak{w}})(T)$. 
 The case is obvious from Theorem \ref{diam4}.
   %  The proof is analogous of case 1 in the proof of Theorem \ref{diam4}.\\
   % { Case 2}: Suppose $v_c \in \Omega_{\mathfrak{w}}(T)$.
   
   % The proof is analogous of  case 2 in the proof of Theorem \ref{diam4}.\\
   % \textcolor{red}{Case 1. and Case 2. proofs, follows from Theorem  \ref{diam4}}\\
   \noindent {Case 2}: Suppose $v_c\in \Omega_{\mathfrak{n}}(T)$.\\
    In this case, $g=w(v_c)=\deg(v_c)g$. Then $(\deg(v_c)-1)g=0$. Therefore $\deg(v_c)-1\equiv 0 \pmod {o(g)}.$ Also, $e(\Gamma)\equiv 0 \pmod {o(g)} $ implies $\gcd(\deg(v_c)-1, e(\Gamma))\equiv 0 \pmod {o(g)}$. Now $o(g)\ne 1$ implies $\gcd(\deg(v_c)-1, e(\Gamma))\ne 1$. Therefore, $T$ satisfies (iii).
    
    Conversely, if $T$ satisfies (i), then by Theorem \ref{internal}, $T$ is $\Gamma$-vertex magic.
    
   Suppose $T$ satisfies (ii). Let $z_i$ be a non pendant vertex of $T$ and $r_i$ denote the number of pendant vertices adjacent to $z_i$. By Lemma \ref{lemma1.6}, label the $r_i$ pendant neighbors of $z_i$ such that whose label sum is zero. Let $g\in \Gamma $ such that $o(g)=e(\Gamma)$. Assign $g$ to all the non-pendant vertices and $(2-\deg(v_c))g$ to the unique pendant neighbor of $v_c$.
    
    Suppose $T$ satisfies (iii). Let $\gcd(\deg(v_c)-1,e(\Gamma))=m>1$. Let $p$ be a prime divisor of $m$. Choose an element $g$ of $\Gamma$ of order $p$.  Hence $\deg(v_c)g=g$. Assign the label $g$  to all non-pendant vertices $z_i\ne v_c$ of $T$, and label an element $g', g'\neq g$ to  $v_c$ and label the pendant neighbors of  $u_i$ in such a way that their labels sum is $g-g'$. This gives an $\Gamma$-vertex magic labeling of $T$.
\end{proof}

\begin{theorem}
    	Let $\Gamma$ be an infinite Abelian torsion-free    group and $T$ be a tree of diameter $4$ with the central vertex $v_c$. Then $T$ is  $\Gamma$-vertex magic if and only if $T$ satisfies one of the criteria below:\\
		\textnormal{(i)} Any non-pendant vertex of $T$ is in $\Omega_{\mathfrak{s}}(T)$.\\
		\textnormal{(ii)} $v_c \in \Omega_{\mathfrak{w}}(T)$  and  all other  non-pendant vertices are in $\Omega_{\mathfrak{s}}(T)$.
\end{theorem}
\begin{proof}
    Let $T$ be $\Gamma$-vertex magic  tree with a labeling $\ell$ and magic constant $g$. There are three cases.\\
    { Case 1}: $v_c \in \Omega_{\mathfrak{s}}(T)$. 
    This case is obvious from Theorem \ref{diam4}.\\
   { Case 2}: $v_c$ is not a support vertex.
   By Theorem \ref{thm2.2}, all the  non-pendant vertices other than $v_c$ are strong support vertices.\\
    { Case 3}: $v_c \in \Omega_{\mathfrak{n}}(T)$.
     In this case, $w(v_c)=\deg(v_c)g=g$ and then, $(\deg(v_c)-1)g=0$. Hence, $\deg(v_c)=1$, which is not possible. Therefore $T$ is not $\Gamma$-vertex magic.
    
Conversely, suppose that $T$ satisfies any of the two conditions.
     
If $T$ satisfies (i), then the result follows from Theorem \ref{internal}.
If $T$ satisfies (ii), then choose  $g \in \Gamma\setminus\{0\}$ and label all the non-pendant support vertices of $T$ by $g$ and label the unique pendant neighbor of $v_c$ by $(2-\deg(v_c))g$.  Label the pendant neighbors of $z_i$ using the $r_i$ number of group elements such that whose sum is 0. Then $T$ is $\Gamma$-vertex magic.
\end{proof}
 \begin{theorem}
    	Let $\Gamma$ be an infinite Abelian  group, which is neither a torsion group nor a torsion-free group and $T$ be a tree of diameter $4$ with the central vertex $v_c$. Then $T$ is  $\Gamma$-vertex magic if and only if $T$ satisfies  one  of the following conditions:\\
		\textnormal{(i)} Any non-pendant vertex of $T$ is in $\Omega_{\mathfrak{s}}(T)$.\\
		\textnormal{(ii)} $v_c \in \Omega_{\mathfrak{w}}(T)$ and  all other non-pendant vertices are in $\Omega_{\mathfrak{s}}(T)$.\\
		\textnormal{(iii)} $v_c$ either in $\Omega_{\mathfrak{n}}(T)$ or in $\Omega_{\mathfrak{p}}(T)$ and $\gcd(\deg(v_c)-1,e(\Gamma))\neq 1.$
\end{theorem}
\begin{proof}
  Let $T$ be $\Gamma$-vertex magic  tree with a labeling $\ell$ and  a magic constant $g$. There are three cases.\\
    { Case 1}: $v_c \in \Omega_{\mathfrak{s}}(T)$. By Theorem \ref{thm2.2}, $T$ satisfies (i).\\
    { Case 2}: $v_c \in \Omega_{\mathfrak{w}}(T)$.
    By Theorem \ref{thm2.2}, $T$ satisfies (ii).\\
    { Case 3}: $v_c \in \Omega_{\mathfrak{n}}(T)$ or $v_c \in \Omega_{\mathfrak{p}}(T)$.\\
In this case, $g=w(v_c)=\deg(v_c)g$. Then $(\deg(v_c)-1)g=0$. Since for any tree of diameter of $4$ and $\deg(v_c)\ge 2$, we have $\deg(v_c)-1\ne 0$. Therefore, $g$ cannot be an element of infinite order. Now, if $g$ is of finite order $\deg(v_c)-1\equiv 0 \pmod {o(g)}$. Also, $e(\Gamma)\equiv 0 \pmod {o(g)} $ implies $\gcd(\deg(v_c)-1, e(\Gamma))\equiv 0 \pmod {o(g)}$. Now, $g \ne 0$ implies $o(g)\ne 1$. Hence $\gcd(\deg(v_c)-1, e(\Gamma))\ne 1$. Therefore, $T$ satisfies (iii).

Conversely, suppose $T$ satisfies any of the three conditions. 

If $T$ satisfies (i), the result follows from Theorem \ref{internal}.  Suppose $T$ satisfies (ii). Choose $g\in \Gamma$ such that $o(g)$ is not finite and label all the  vertices in or $v_c \in \Omega_{\mathfrak{n}}(T)$ by $g$ and label the unique pendant neighbor $v_c$ by $(2-\deg(v_c))g$.  Label the pendant neighbors of $z_i$ by the $r_i$ number of group elements such that their label sum is $0$. Then $T$ is $\Gamma$-vertex magic. Suppose $T$ satisfies (iii). Let $\gcd(d(v_c)-1,e(\Gamma))=m>1$. Let $p$ be a prime divisor of $m$. Choose $g\in \Gamma$ such that $o(g)=p$. Hence $d(v_c)g=g$. Assign $g$ to all the non-pendant vertices $u_i\ne v_c$, assign an element $g'\in \Gamma-\{0,g\}$ to $v_c$ and label the pendant neighbors of  $u_i$ such that their label sum is $g-g'$. This gives an $\Gamma$-vertex magic labeling of $T$.
\end{proof}

Thus, we characterized all $\Gamma$-vertex magic trees of diameter $4$ for infinite group with finite torsion elements. Now, we proceed to characterize all $\Gamma$-vertex magic trees of diameter $5$. 
\begin{figure}[htbp]
\centering
   \begin{tikzpicture}[scale=0.45]
   \filldraw[lightgray,line width=5pt,rounded corners=4pt]  (0,0.1) ellipse (1cm and 0.5cm);
    \filldraw[lightgray,line width=5pt,rounded corners=4pt]  (0,-4) ellipse (4.5cm and 0.5cm);
     \filldraw[lightgray,line width=5pt,rounded corners=4pt]  (-4.5,-8) ellipse (1.8cm and 0.5cm);
     \filldraw[lightgray,line width=5pt,rounded corners=4pt]  (-0.5,-8) ellipse (1.8cm and 0.5cm);
     \filldraw[lightgray,line width=5pt,rounded corners=4pt]  (4.5,-8) ellipse (1.8cm and 0.5cm);
    \filldraw[fill=black](0,0) circle (0.2cm);
    \filldraw[fill=black](-4,-4) circle (0.2cm);
    \filldraw[fill=black](-1,-4) circle (0.2cm);
    \filldraw[fill=black](4,-4) circle (0.2cm);
    \filldraw[fill=black](-6,-8) circle (0.2cm);
     \filldraw[fill=black](-5,-8) circle (0.2cm);
      \filldraw[fill=black](-3,-8) circle (0.2cm);
       \filldraw[fill=black](-2,-8) circle (0.2cm);
        \filldraw[fill=black](-1,-8) circle (0.2cm);
         \filldraw[fill=black](1,-8) circle (0.2cm);
          \filldraw[fill=black](6,-8) circle (0.2cm);
           \filldraw[fill=black](3,-8) circle (0.2cm);
            \filldraw[fill=black](4,-8) circle (0.2cm);
            \draw(0,0)--(4,-4);
\draw(0,0)--(-4,-4);
\draw(0,0)--(-1,-4);
\draw(-6,-8)--(-4,-4);
 \draw(-5,-8)--(-4,-4); 
 \draw(-3,-8)--(-4,-4);
 \draw(-2,-8)--(-1,-4);
 \draw(-1,-8)--(-1,-4);
 \draw(1,-8)--(-1,-4);
 \draw(6,-8)--(4,-4);
  \draw(3,-8)--(4,-4);
   \draw(4,-8)--(4,-4);
 \filldraw[fill=black](0.5,-4) circle (0.02cm); 
  \filldraw[fill=black](1.5,-4) circle (0.02cm); 
   \filldraw[fill=black](2.5,-4) circle (0.02cm); 
\node at (-4,-8){$\cdots$};
 \node at (0,-8){$\cdots$}; 
 \node at (5,-8){$\cdots$};
 \node at (0,0.6){$v_{c_1}$};

  \filldraw[lightgray,line width=5pt,rounded corners=4pt]  (14,0.1) ellipse (1cm and 0.5cm);
    \filldraw[lightgray,line width=5pt,rounded corners=4pt]  (14,-4) ellipse (4.5cm and 0.5cm);
     \filldraw[lightgray,line width=5pt,rounded corners=4pt]  (9.5,-8) ellipse (1.8cm and 0.5cm);
     \filldraw[lightgray,line width=5pt,rounded corners=4pt]  (13.5,-8) ellipse (1.8cm and 0.5cm);
     \filldraw[lightgray,line width=5pt,rounded corners=4pt]  (18.5,-8) ellipse (1.8cm and 0.5cm);
    \filldraw[fill=black](14,0) circle (0.2cm);
    \filldraw[fill=black](10,-4) circle (0.2cm);
    \filldraw[fill=black](13,-4) circle (0.2cm);
    \filldraw[fill=black](18,-4) circle (0.2cm);
    \filldraw[fill=black](8,-8) circle (0.2cm);
     \filldraw[fill=black](9,-8) circle (0.2cm);
      \filldraw[fill=black](11,-8) circle (0.2cm);
       \filldraw[fill=black](12,-8) circle (0.2cm);
        \filldraw[fill=black](13,-8) circle (0.2cm);
         \filldraw[fill=black](15,-8) circle (0.2cm);
          \filldraw[fill=black](20,-8) circle (0.2cm);
           \filldraw[fill=black](17,-8) circle (0.2cm);
            \filldraw[fill=black](18,-8) circle (0.2cm);
            \draw(14,0)--(18,-4);
\draw(14,0)--(10,-4);
\draw(14,0)--(13,-4);
\draw(8,-8)--(10,-4);
 \draw(9,-8)--(10,-4); 
 \draw(11,-8)--(10,-4);
 \draw(12,-8)--(13,-4);
 \draw(13,-8)--(13,-4);
 \draw(15,-8)--(13,-4);
 \draw(20,-8)--(18,-4);
  \draw(17,-8)--(18,-4);
   \draw(18,-8)--(18,-4);
 \filldraw[fill=black](14.5,-4) circle (0.02cm); 
  \filldraw[fill=black](15.5,-4) circle (0.02cm); 
   \filldraw[fill=black](16.5,-4) circle (0.02cm); 
\node at (10,-8){$\cdots$};
 \node at (14,-8){$\cdots$}; 
 \node at (19,-8){$\cdots$};
 \node at (14,0.6){$v_{c_2}$};
 \draw(0,0)--(14,0);
  \node at (-2,0){$\Omega_{\mathfrak{n}}$};
 \node at (-7,-4){$\Omega_{\mathfrak{s}}\cup \Omega_{\mathfrak{w}}$};
  \node at (5,-9.2){$\Omega_{\mathfrak{p}}$};
  \node at (-0.1,-9.2){$\Omega_{\mathfrak{p}}$};
  \node at (-4,-9.2){$\Omega_{\mathfrak{p}}$};
   \node at (16,0){$\Omega_{\mathfrak{n}}$};
 \node at (21,-4){$\Omega_{\mathfrak{s}}\cup \Omega_{\mathfrak{w}}$};
  \node at (19,-9.2){$\Omega_{\mathfrak{p}}$};
  \node at (13.9,-9.2){$\Omega_{\mathfrak{p}}$};
  \node at (10,-9.2){$\Omega_{\mathfrak{p}}$};
\end{tikzpicture}
\caption{A typical tree with diameter $5$}
\end{figure}
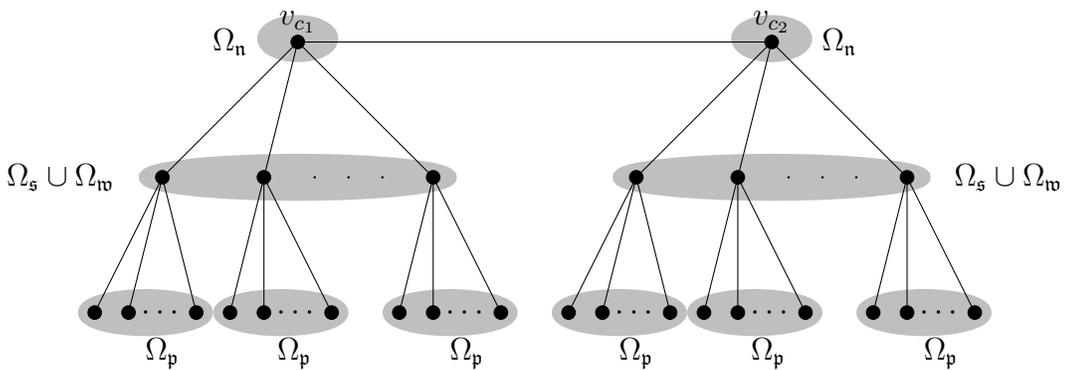
\begin{theorem}
Let $\Gamma$ be an  infinite Abelian group, which is not a torsion group and $T$ be a tree of diameter $5$ such that both central vertices $v_{c_1}$ and $v_{c_2}$ are not support vertices. Then $T$ is $\Gamma$-vertex magic if and only if
$\deg(v_{c_i})\ne 2$. 
\end{theorem}
\begin{proof}
	Suppose $T$ is  $\Gamma$-vertex magic with a labeling $\ell$ and  magic constant $g$. Clearly $g\ne0$ and $\ell(v)=g$ for all the support vertices $v$. Then $g=w(v_{c_{1}})= (\deg(v_{c_1})-1)g+\ell(v_{c_{2}})$. Hence  $\ell (v_{c_{2}})=(2-\deg(v_{c_{1}}))g$. Since $\ell(v_{c_{2}})\ne 0$, we have $~(2-\deg(v_{c_1}))g\ne 0$. Therefore  $\deg(v_{c_1})\ne 2$.
	 Similarly, we get  $\deg(v_{c_{2}})\ne 2$.
	 
	 Conversely, if $\deg(v_{c_i})\ne 2$, then choose a non-torsion element $g$ of $\Gamma$. Let $g'=(\deg(v_{c_1})-1)g$ and $g''=(\deg(v_{c_2})-1)g$. Since $\deg(v_{c_i})\ne 2$, we have $g-g'$ and $g-g''$ are non-zero. Now define $\ell(v_{c_1})=g-g''$, $\ell(v_{c_2})=g-g',$ and $\ell(v)=g$, for any support vertex $v$ of $T$.  If $u$ is a support vertex adjacent to $v_{c_1}$, label the pendant neighbors of $u$ such that their label sum is $g''$. If $u$ is a support vertex adjacent to $v_{c_2}$, label the pendant neighbors of $u$ such that their label sum is $g'$. Thus, $T$ is  $\Gamma$-vertex magic.
	\end{proof}

 \begin{theorem}
Let $\Gamma$ be an infinite Abelian group, which is not a torsion group and $T$ be a tree of diameter $5$ such that both its central vertices $v_{c_1}$ and $v_{c_2}$ are in $\Omega_{\mathfrak{s}}\cup \Omega_{\mathfrak{w}}$. Then $T$ is $\Gamma$-vertex magic if and only if  $T$ has no vertex of degree $2$.
\end{theorem}
\begin{proof}
Suppose $T$ is $\Gamma$-vertex magic with a labeling $\ell$ and magic constant $g$. Since $v_{c_1}$ and $v_{c_2}$  are support vertices, $\ell(v_{c_1})=g=\ell(v_{c_2})$, and both $v_{c_1}$ and $v_{c_2}$ have degrees more than 2.

Now if $T$ has a vertex $y$ with $\deg(y)=2$, then $y$ is a support vertex and hence $N(y)=\{v_{c_i},z\}$, where $z$ is a pendant vertex. Therefore,  $g=w(y)=\ell(v_{c_i})+\ell(z)=g+\ell(z)$ and thus $\ell(z)=0$,  a contradiction. Hence $T$ has no vertex of degree 2. 	

\par Conversely, suppose $T$ has no vertex of degree $2$. Since $\Gamma$ is not torsion group, let $g$ be any non-torsion element in $\Gamma$.   Let $r_i$ denote the number of pendant neighbors adjacent to $v_{c_i}$ and $g'=(\deg(v_{c_1})-r_1)g$ and $g''=(\deg(v_{c_2})-r_2)g$. Label all the support vertices of $T$ by $g$. Label the pendant neighbors of $v_{c_1}$ such that their label sum is $g-g'$. Similarly, label the pendant neighbors of $v_{c_2}$ such that their label sum is $g-g''$. If $y$ is a support vertex adjacent to $v_{c_i}$, then label all the pendant neighbors of $y$ such that their label sum is $0$.
\end{proof}

\begin{theorem}
		Let $\Gamma$ be an infinite torsion-free Abelian group and $T$ be a tree of diameter $5$ such that  $v_{c_{1}}\in \Omega_{\mathfrak{s}}\cup \Omega_{\mathfrak{w}}$ and  $v_{c_{2}}\not\in \Omega_{\mathfrak{s}}\cup \Omega_{\mathfrak{w}}$. Then $T$ is not $\Gamma$-vertex magic.
	\end{theorem}
\begin{proof}
Let $\ell$ be an $\Gamma$-vertex magic labeling with a magic constant $g$. If $v\ne v_{c_2}$ with $\deg(v)\ge 2$ is a support vertex, then $\ell(v)=g$. 
Now, $g=w(v_{c_2})=\deg(v_{c_2})g$ and then, $(\deg(v_{c_2})-1)g=0.$ But $\deg(v_{c_2})$ cannot be equal to 1. Thus, $o(g)$ cannot be infinite.
\par Conversely, if $T$ is tree such that $v_{c_1}$ is not adjacent to a weak support vertex, then choose a non-torsion element $g\in \Gamma$. 
		Label $\ell(v)=g$, for all support vertices $v$ of $T$.  Let $r$ be the number of pendants adjacent to $v_{c_1}$ and $g' = (\deg(v_{c_{1}})-r)g$. 
	
	%  and choose $ h\in \Gamma- %\{0,g,g'\}$.\\
	%	If $v_{c_1}$ is a strong support vertex, 
		Then define $\ell(v_{c_2})=g'$ and label the pendant neighbors of $v_{c_1}$ in such a way that their label sum is $2g-2g'$.
		
	%	If $v_{c_1}$ is a weak support vertex with unique pendant neighbor $z$, then choose $h'\in \Gamma-\{0,g,g',h\}$. Now, define $\ell(v_{c_2})=h'$ and $\ell(z)=2g-(g'+h')$. Then $w(v_{c_1})=w(c_{v_2})=g$.
		
	Finally, if $y$ is a support vertex adjacent to $v_{c_1}$, then label pendant neighbors of $y$ such that their label sum is $0$. If $z$ be a support vertex adjacent to $v_{c_2}$ then label the pendant neighbors of $z$ such that their label sum is $g-g'$.  This gives an $\Gamma$-vertex magic labeling of $T$. 
	\end{proof}

\begin{theorem}	Let $\Gamma$ be an infinite Abelian group, which is not a torsion-free group and $T$ be a tree of diameter $5$ such that   $v_{c_{1}}\in \Omega_{\mathfrak{s}}(T)\cup \Omega_{\mathfrak{w}}(T)$ and  $v_{c_{2}}\not\in \Omega_{\mathfrak{s}}(T)\cup \Omega_{\mathfrak{w}}(T)$. Then $T$ is $\Gamma$-vertex magic if and only if $T$ satisfies the following conditions:\\
\textnormal{(i)} $v_{c_1}$ is not adjacent to any vertices in  $ \Omega_{\mathfrak{w}}(T)$,\\
\textnormal{(ii)} $\gcd(\deg(v_{c_2})-1,e(\Gamma))\ne 1$. 
\end{theorem}
\begin{proof}
Suppose $T$ is  $\Gamma$-vertex magic  with a labeling $\ell$ and  magic constant $g$. Any vertex $y\ne v_{c_2}$ with $\deg(y)\ge 2$ is a support vertex and hence $\ell(y)=g$. 

Let $v_{c_1}$ is adjacent to $u$, $u \in  \Omega_{\mathfrak{w}}(T)$ and if $z$ is the unique  pendant vertex adjacent to $y$, then $g=w(y)=\ell(z)+\ell(v_{c_1})=\ell(z)+g.$ Hence $\ell(z)=0$, which is a contradiction. Therefore, $T$ satisfies (i).
		 
Now, $g=w(v_{c_2})=\deg(v_{c_2})g$ and then, $(\deg(v_{c_2})-1)g=0.$ But $\deg(v_{c_2})$ cannot be equal to 1. Thus, $o(g)$ cannot be infinite. Therefore, $g$  is of finite order and $o(g)$ must divide $\deg(v_{c_2})-1$. On the other hand, $o(g)$ divides $e(\Gamma)$ and hence $o(g)$ divides $\gcd(\deg(v_{c_2})-1,e(\Gamma))$. Since $g\ne 0, o(g)\ne 1$ and hence $\gcd(\deg(v_{c_2})-1,e(\Gamma))\ne 1$. Therefore, $T$ satisfies (ii).

		\par Conversely, suppose $T$ satisfies (i) and (ii). 
		Let $\gcd(\deg(v_c)-1,e(\Gamma))=m$ and   $p$ be a prime divisor of $m$.  Let  $g\in \Gamma$ such that $o(g)=p$. Label $\ell(v)=g$ for all support vertices $v$. Let $g' = (\deg(v_{c_{1}})-1)g$  and choose $ h\in \Gamma- \{0,g,g'\}$. If $v_{c_1}$ is a strong support vertex, then define $\ell(v_{c_1})=h$ and label the pendant neighbors of $v_{c_1}$ in such a way that their label sum is $g-(g'+h)$. If $v_{c_1}$ is a weak support vertex with a pendant neighbor $z$, choose $h'\in \Gamma-\{0,g,g',h\}$. Define $\ell(v_{c_2})=h'$ and $\ell(z)=g-(g'+h')$. Then $w(v_{c_1})=w(c_{v_2})=g$.
		
	Now, label the remaining pendant vertices in such a way that $w(z)=g$ for all the remaining support vertices, by Lemma \ref{lemma1.6}. This gives an $\Gamma$-vertex magic labeling of $T$. 
\end{proof}

Therefore, above theorems provide a complete characterization of all $\Gamma$- vertex magic trees of diameter 5 of infinite group with finite torsion element.

 \section{Conclusion}
In this article,  the group distance magicness of various product graphs are proved and then a characterization of $\Gamma$-vertex magicness of trees of diameter at most $5$, is obtained, where  $\Gamma$ is any infinite Abelian group with finitely many torsion elements. However the following problem is still open. 
\vskip 0.1cm \noindent 
{\bf Problem 1.} If $\Gamma$ is an infinite Abelian group with infinitely many torsion elements, characterize $\Gamma$-vertex magicness of trees of diameter at most $5$.

% \nocite{*}
% \bibliographystyle{plain}
% \bibliography{karthik.bib}

\begin{thebibliography}{9}
\bibitem{MR3445306} G. Chartrand, L. Lesniak, and P. Zhang, \textit{Graphs \& Digraphs}, Textbooks in Mathematics, CRC Press, Boca Raton, FL, 2016.
\bibitem{MR2817074} R. Hammack, W. Imrich, and S. Kla\u{v}zar, \textit{Handbook of Product Graphs}. Discrete Mathematics and its Applications, CRC Press, Boca Raton, FL, 2011.
\bibitem{MR0356988} I. N. Herstein, \textit{Topics in Algebra}, Xerox College Publishing, Lexington, Mass.-Toronto, 1975.
\bibitem{MR3013201} A. M. Marr and W. D. Wallis, \textit{Magic Graphs}, Springer, New York, 2013.
\bibitem{MR281659} R. Frucht and F. Harary,  On the corona of two graphs, \textit{Aequationes Math.,} \textbf{4}:322–325, 1970.
\bibitem{MR1853016} S Lee, I Wen, and H Sun, On group-magic graphs, \textit{J. Combin. Math. Combin. Comput.,} \textbf{38}:197–207, 2001.
\bibitem{MR4145437} N. Kamatchi, K. Paramasivam, A. V. Prajeesh, K. M Sabeel, and S. Arumugam, On group vertex magic graphs, \textit{AKCE Int. J. Graphs Comb.,} \textbf{17(1)}:461–465, 2020.
\bibitem{MR4517978} S. Balamoorthy, S. V. Bharanedhar, and N. Kamatchi, On the products of group vertex magic graphs, \textit{AKCE Int. J. Graphs Comb.,} \textbf{19(3)}:268–275, 2022.
\bibitem{MR4515472} K. M. Sabeel, K. Paramasivam, A.V. Prajeesh, N. Kamatchi, and S. Arumugam, A characterization of group vertex magic trees of diameter up to 5, \textit{Australas. J. Combin.,} \textbf{85}:49--60, 2023.
\bibitem{MR1668059} J. A. Gallian, A dynamic survey of graph labeling, \textit{Electron. J. Combin.,} \#DS6, 2022.
%\bibitem{MR4179413} D.Andrica, S. R\u{a}dulescu, and G.C. \c{T}urca\c{s}, \textit{The exponent of a group: properties, computations and applications}, Discrete Mathematics and Applications, Springer Optim. Appl., \textbf{165}, 57-108, 2020.
\bibitem{1} K. M Sabeel, A.V. Prajeesh, and K. Paramasivam, A characterization for $V_4$-vertex magicness of trees with diameter 5,  \textit{In Computational Sciences
- Modelling, Computing and Soft Computing}, 243–249, Springer, 2021.
\end{thebibliography}

\end{document}